\documentclass[a4paper, 12pt]{amsart}
\usepackage{cite}
\usepackage{graphicx}
\usepackage{color}
\usepackage[T1]{fontenc}

\usepackage{amsmath}
\usepackage{amssymb}
\usepackage{amsthm}
\usepackage{bm}
\usepackage{bbm}
\usepackage{eucal}
\usepackage{braket}
\usepackage{stmaryrd}
\usepackage{mathrsfs}
\usepackage{mathtools}

\usepackage[all,2cell]{xy}
\UseAllTwocells
\usepackage{youngtab}
\usepackage{tikz}
\usetikzlibrary{intersections,calc,arrows.meta}

\usepackage[hang,small,bf]{caption}
\usepackage[subrefformat=parens]{subcaption}
\calclayout

\theoremstyle{plain}
\newtheorem{thm}{Theorem}[section]
\newtheorem{lem}[thm]{Lemma}
\newtheorem{prop}[thm]{Proposition}

\theoremstyle{definition}
\newtheorem{defn}[thm]{Definition}%[section]

\theoremstyle{remark}
\newtheorem{rem}[thm]{Remark}%[section]
\newtheorem{exam}[thm]{Example}
\makeatletter
    
    \@addtoreset{equation}{section}
\makeatother

\usepackage{xr-hyper}
\usepackage{hyperref}  %comment out in setup.tex

\makeatletter

\newcommand*{\addFileDependency}[1]{% argument=file name and extension
\typeout{(#1)}% latexmk will find this if $recorder=0
\@addtofilelist{#1}
\IfFileExists{#1}{}{\typeout{No file #1.}}
}\makeatother

\newcommand{\supp}{\mathrm{supp}}

\newcommand{\rmRe}{\mathrm{Re}}
\newcommand{\rmIm}{\mathrm{Im}}
\newcommand{\dee}{\partial}

\newcommand{\arcsinh}{\mathrm{arcsinh}}

\newcommand{\what}[1]{\widehat{#1}}
\newcommand{\wtilde}[1]{\widetilde{#1}}

\newcommand{\ol}[1]{\overline{#1}}

\newcommand{\Cov}{\mathrm{Cov}}
\newcommand{\Var}{\mathrm{Var}}

\newcommand{\law}{\mathrm{(law)}}
\newcommand{\leb}{\mathsf{m}}
\newcommand{\weakConv}{\xrightarrow{\mathrm{weak}}}

\newcommand{\gffD}{\Gamma}
\newcommand{\harm}{\mathfrak{u}}
\newcommand{\harmAt}[1]{\harm_{#1}}
\newcommand{\harmMAt}[2]{\harm^{[#1]}_{#2}} % M for multiple
\newcommand{\cpxHarmMAt}[2]{\wtilde{\harm}^{[#1]}_{#2}}
\newcommand{\gffB}{h}
\newcommand{\harmCl}{\mathsf{u}}
\newcommand{\harmClMAt}[2]{\harmCl^{[#1]}_{#2}}
\newcommand{\harmHydro}{\harm^{[\infty]}}
\newcommand{\harmHydroAt}[1]{\harmHydro_{#1}}

\newcommand{\SLEg}{g}
\newcommand{\SLEgAt}[1]{\SLEg_{#1}}
\newcommand{\mSLEg}[1]{\SLEg^{[#1]}}
\newcommand{\mSLEgAt}[2]{\mSLEg{#1}_{#2}}

\newcommand{\SLEDom}{\bD}
\newcommand{\SLEDomAt}[1]{\SLEDom_{#1}}
\newcommand{\mSLEDom}[1]{\SLEDom^{[#1]}}
\newcommand{\mSLEDomAt}[2]{\mSLEDom{#1}_{#2}}

\newcommand{\SLEHull}{K}
\newcommand{\SLEHullAt}[1]{\SLEHull_{#1}}
\newcommand{\mSLEHull}[1]{\SLEHull^{[#1]}}
\newcommand{\mSLEHullAt}[2]{\mSLEHull{#1}_{#2}}

\newcommand{\empMeas}{\Xi}
\newcommand{\DysonEmpAtOf}[3]{\empMeas^{[#1]}_{#2}(#3)}
\newcommand{\DysonEmpAt}[2]{\DysonEmpAtOf{#1}{#2}{\cdot}}

\newcommand{\diracMeasAtOf}[2]{\delta_{#1}(#2)}
\newcommand{\diracMeasAt}[1]{\diracMeasAtOf{#1}{\cdot}}

\newcommand{\cons}{\mathfrak{w}}
\newcommand{\consAt}[1]{\cons_{#1}}

\newcommand{\bC}{\mathbb{C}}
\newcommand{\bD}{\mathbb{D}}
\newcommand{\bE}{\mathbb{E}}

\newcommand{\bH}{\mathbb{H}}

\newcommand{\bN}{\mathbb{N}}

\newcommand{\bP}{\mathbb{P}}

\newcommand{\bR}{\mathbb{R}}
\newcommand{\bS}{\mathbb{S}}

\newcommand{\cF}{\mathcal{F}}

\newcommand{\cH}{\mathcal{H}}

\newcommand{\sfM}{\mathsf{M}}

\newcommand{\sfU}{\mathsf{U}}

\newcommand{\bbmi}{\mathbbm{i}}

\def\D{\mathbb{D}}

\def\N{\mathbb{N}}

\def\C{\mathbb{C}}
\def\R{\mathbb{R}}

\def\1{{\bf 1}}

\def\cO{\mathcal{O}}

\title[Radial multiple SLE, GFF, and hydrodynamic limit]{Coupling of radial multiple SLE with Gaussian free field, and the hydrodynamic limit}
\date{\today}
\author{Makoto Katori}
\address{Department of Physics, Faculty of Science and Engineering, Chuo University, Kasuga, Bunkyo-ku, Tokyo 112-8551, Japan}
\email{makoto.katori.mathphys@gmail.com}

\author{Shinji Koshida}
\address{Department of Mathematics and Systems Analysis, Aalto University, Espoo, Finland}
\email{shinji.koshida@aalto.fi}

\author{Chizuru Soukejima}
\address{Department of Physics, Faculty of Science and Engineering, Chuo University, Kasuga, Bunkyo-ku, Tokyo 112-8551, Japan}
\email{a22.55ak@g.chuo-u.ac.jp}

\author{Raian Suzuki}
\address{Department of Physics, Faculty of Science and Engineering, Chuo University, Kasuga, Bunkyo-ku, Tokyo 112-8551, Japan}
\email{a22.sjce@g.chuo-u.ac.jp}

\begin{document}

\begin{abstract}
Schramm--Loewner evolution (SLE) has been one of the central topics in the probabilistic study of two-dimensional critical systems.
It is a random curve in two dimensions to which a cluster interface in a critical lattice system is supposed, or has been proved, to converge.
The most archetypical setting for SLE is called chordal, where a random curve evolves in a simply-connected domain from a boundary point to another, whereas in its variant called radial, a random curve evolves from a boundary point to a distinguished interior point.
Multiple SLE is a variant to another direction; it deals with multiple random curves, and it is a natural direction as there are certainly multiple cluster interfaces found in critical lattice systems.
In this paper, we study multiple SLE in the radial setting, namely, radial multiple SLE.

We report two main results. One is regarding the local set coupling between radial multiple SLE and Gaussian free field (GFF).
This sort of coupling between SLE and GFF has been extensively studied in the chordal setting, and serves as a foundation for many recent developments.
We show that the coupling between radial multiple SLE and GFF occurs if and only if the radial multiple SLE is driven by the circular Dyson Brownian motions.

The circular Dyson Brownian motions are a typical example of stochastic log-gases.
This fact motivates us to study the hydrodynamic limit of the corresponding radial multiple SLE, which refers to the dynamical law of large numbers, when the number of curves tends to infinity.
In our other main result, we provide explicit description of the hydrodynamic limit.

\end{abstract}

\maketitle

\tableofcontents

\section{Introduction}
\subsection{Background}
\subsubsection*{Schramm--Loewner evolution}
Schramm--Loewner evolution (SLE)~\cite{Schramm2000} is a random continuous curve in a two-dimensional domain that is conformally invariant and satisfies the domain Markov property.
In the most studied {\it chordal} case, the curve evolves between two boundary points.
Due to the conformal invariance, we may assume that the domain is the complex upper-half plane $\bH=\{z\in \bC| \rmIm z>0\}$, and that the curve evolves from $0$ to $\infty$.
Thanks to the Loewner theory from complex analysis, such a curve is governed by the chordal Loewner equation 
\begin{align}
\label{eq:chord_loewner}
	\frac{d}{dt}\SLEgAt{t}(z) &= \frac{2}{\SLEgAt{t}(z)-X_{t}},\quad t\geq 0, \\
	\SLEgAt{0}(z) &= z\in \bH, \notag
\end{align}
where $(X_{t}:t\geq 0)$ is a continuous stochastic process on $\bR$ called the driving process.
The solution $(\SLEgAt{t}:t\geq 0)$ is called a Loewner chain and, at each $t\geq 0$, 
\begin{align*}
	\SLEgAt{t}\colon \bH_{t} \to \bH
\end{align*}
is conformal\footnote{
In this paper, we call a bijective conformal map simply a conformal map.
}, where $\bH_{t}\subset \bH$ is the set of $z\in \bH$ for which (\ref{eq:chord_loewner})
has a solution until $t\geq 0$.
The original curve, which is denoted by 
$\eta\colon [0,\infty) \to \ol{\bH}$,
is recovered from the Loewner chain as 
\begin{align*}
	\eta_{t} = \lim_{\epsilon \downarrow 0}\SLEgAt{t}^{-1}(X_{t} 
	+\bbmi \epsilon),\quad t\geq 0,
\end{align*}
where $\bbmi=\sqrt{-1}$.

Remarkably, the driving process must be $X_{t}=\sqrt{\kappa}B_{t}$, $t\geq 0$
where $(B_{t}:t\geq 0)$ is a standard Brownian motion and $\kappa\geq 0$ when the law of $\eta$ is conformally invariant and satisfies the domain Markov property~\cite{Schramm2000}.
With this choice of a driving process, the curve $\eta$ is indeed continuous~\cite{RohdeSchramm2005}.
In other words, there is a unique SLE for each $\kappa \geq 0$.
More detailed and thorough account of SLE can be found in~\cite{Lawler2005,Katori2015,kemppainen2017schramm}.

\subsubsection*{Critical lattice systems and conformal field theory}
SLE was introduced in the study of critical lattice models in two dimensions 
as a candidate for a scaling limit of a domain interface.
In fact, the two characterizing properties of SLE, conformal invariance and domain Markov property, are naturally expected for a scaling limit of a domain interface.
There are by now various critical lattice models whose domain interface is proved to converge to an SLE~\cite{Smirnov2001,LawlerSchrammWerner2004,schramm2005harmonic,SchrammSheffield2009,ChelkakDuminil-CopinHonglerKemppainenSmirnov2014}.

A critical lattice system in two dimensions is also conjectured to give a conformal field theory (CFT)~\cite{BelavinPolyakovZamolodchikov1984,belavin1984infinite} in a scaling limit.
Therefore, one might expect an interrelation between SLE and CFT as both capture some properties of scaling limits of critical lattice systems.
Such an interrelation was studied right after the introduction of SLE and dubbed the SLE/CFT-correspondence~\cite{BauerBernard2002,BauerBernard2003,BauerBernard2004, FriedrichWerner2003, FriedrichKalkkinen2004,KontsevichSuhov2007, Kytola2007,Dubedat2015a,Dubedat2015b}.
It would be more appropriate to call it a {\it coupling} rather than {\it correspondence}, though, because of the reason explained next.

\subsubsection*{Gaussian free field}
While the above mentioned studies of the SLE/CFT-correspondence were more or less based on the representation theory of the Virasoro algebra as a formulation of CFT,
Schramm and Sheffield~\cite{SchrammSheffield2013} found a purely probabilistic instance of the SLE/CFT-correspondence using Gaussian free field (GFF), which is a probabilistic avatar of the CFT of massless free bosons.
In their work, the SLE/CFT-correspondence was realized as a coupling between SLE and GFF.
This coupling between SLE and GFF, nowadays commonly referred to as the local set coupling, was studied further~\cite{Dubedat2009} and serves as a foundation of the Liouville quantum gravity~\cite{DuplantierMillerSheffield2014,Sheffield2016,ang2023conformal,ang2023conformalb}, imaginary geometry~\cite{MillerSheffield2016a, MillerSheffield2016b, MillerSheffield2016c, MillerSheffield2017}, and various other developments~\cite{KangMakarov2013,Koshida2019,KatoriKoshida2020a,KatoriKoshida2020b,ang2022sle,byun2023conformal,alberts2024conformal,aru2024sle,ang2024sle,ang2024integrability}.

\subsubsection*{Multiple SLE}
From the original motivation for SLE, it is natural to extend it to deal with multiple random curves in a two-dimensional domain because there could be certainly found multiple domain interfaces in a critical lattice model.
Such an extension is called multiple SLE.
Similarly to the single curve case, convergence to multiple SLE has been proven for various critical lattice models~\cite{izyurov2015smirnov, izyurov2017critical, BeffaraPeltolaWu2018, Karrila2019, PeltolaWu2019, Karrila2020, izyurov2022multiple, peltola2023crossing}.

A global definition of multiple SLE~\cite{KozdronLawler2007, BeffaraPeltolaWu2018} is concise and conceptual, but its local model based on the Loewner theory is also useful for some purposes.
There are actually two local models of multiple SLE; one is a commuting family of Loewner chains~\cite{Dubedat2007} and the other is a single Loewner chain driven by a multi-particle stochastic process~\cite{BauerBernardKytola2005,Graham2007,Schleissinger2012,RothSchleissinger2017}.

The two approaches are equivalent at least as long as the curves are apart from each other, but the latter is closer to our point of view.
For $N\geq 1$, a multiple SLE with $N$-curves is the Loewner chain that solves the multiple Loewner equation
\begin{align}
\label{eq:chord_multiple_loewner}
	\frac{d}{dt}\SLEgAt{t}(z) &= \sum_{i=1}^{N}\frac{2}{\SLEgAt{t}(z) - X^{(i)}_{t}},\quad t \geq 0, \\
	\SLEgAt{0}(z) &= z\in \bH, \notag
\end{align}
where $(X^{(i)}_{t}:t\geq 0)$, $i=1,\dots, N$ are continuous stochastic processes whose laws are determined by what is called a partition function. See~\cite{BauerBernardKytola2005} for the stochastic differential equations for the driving processes in terms of a partition function.
Equation~(\ref{eq:chord_multiple_loewner}) is obviously an extension of (\ref{eq:chord_loewner}), and similarly, we may find $\bH_{t}\subset\bH$ at each $t\geq 0$ so that 
\begin{align*}
	\SLEgAt{t}\colon \bH_{t} \to \bH
\end{align*}
is conformal.
Furthermore, $N$ curves $\eta^{(i)}\colon [0,\infty) \to \ol{\bH}$, $i=1,\dots, N$ are found as 
\begin{align*}
	\eta^{(i)}_{t} = \lim_{\epsilon \downarrow 0}\SLEgAt{t}^{-1}(X^{(i)}_{t}+\bbmi \epsilon),\quad t\geq 0,\quad i=1,\dots, N.
\end{align*}

Within the framework of multiple SLE, we can at best say that a partition function must satisfy the set of Belavin--Polyakov--Zamolodchikov equations for some consistency~\cite{Dubedat2007}.
However, there is no intrinsic principle that allows us to fix a partition function uniquely. In other words, the laws of the driving processes $(X^{(i)}_{t}:t\geq 0)$, $i=1,\dots, N$ in (\ref{eq:chord_multiple_loewner}) are another input than (\ref{eq:chord_multiple_loewner}) itself.

\subsubsection*{Dyson Brownian motions}
On the other hand, there should be a unique law of multiple curves when we know that they are a scaling limit of domain interfaces of a critical lattice model.
In other words, one should be able to fix a partition function, or equivalently, driving processes by assuming coupling with a specific CFT under specific boundary conditions.
In the previous work by the first two authors~\cite{KK_three_phases_2021}, we studied the local set coupling between multiple SLE and GFF. 
In upshot, the coupling is possible if and only if the driving processes are the Dyson Brownian motions~\cite{Dyson1962} with the parameter determined by the boundary conditions of the GFF.
The Dyson Brownian motions originally arose as a dynamical version of the Gaussian ensemble of random matrices and were extended as stochastic particle systems with logarithmic interaction potential, as known as stochastic log-gases.
Therefore, our result is to intertwine stochastic log-gases including dynamical random matrix theory and multiple SLE through GFF.

\subsubsection*{Hydrodynamic limit}
One of the typical problems regarding a stochastic particle system is its \textit{hydrodynamic limit}.
This is a dynamical analogue of the law of large numbers.
In the random matrix terminology, we can think of it as a large $N$-limit, or a dynamical analogue of Wigner's semi-circular law.
In the case of Dyson Brownian motions, the hydrodynamic limit is well-established~\cite{AndersonGuionnetZeitouni2010}.
The hydrodynamic limit of the corresponding multiple SLE was studied in~\cite{delMonacoSchleissinger2016, delMonacoHottaSchleissinger2018} and exactly solved by Hotta and the first author~\cite{HottaKatori2018}.

\subsection{Radial SLE}
All the above results are regarding the chordal setting.
A variant of the chordal SLE is the {\it radial} SLE, where a curve runs between a boundary point and a distinguished interior point.

Due to the conformal invariance, we may pick the unit disk $\bD=\{z\in \bC||z|<1\}$
as a domain of interest, and suppose that the curve evolves from $1$ to $0$.
Similarly to the chordal case, we rely on the Loewner theory to describe the curve.
The radial Loewner equation reads
\begin{align}
\label{eq:radial_loewner}
    \frac{d}{dt}\SLEgAt{t}(z)&=-\SLEgAt{t}(z)\frac{\SLEgAt{t}(z)+X_{t}}{\SLEgAt{t}(z)-X_{t}},\quad t\geq 0,\\
    \SLEgAt{0}(z)&=z \in \bD, \notag
\end{align}
where the driving process $(X_{t}:t\geq 0)$ is continuous on the unit circle $\bS=\{z\in \bC||z|=1\}$.
Again, for each $t\geq 0$, we define $\SLEDomAt{t}$ as the set of $z\in \bD$ such that the solution of (\ref{eq:radial_loewner}) exists until the time $t$. 
Then, $\SLEgAt{t}\colon \SLEDomAt{t} \to \bD$ is conformal at each $t\geq 0$.
A curve $\eta\colon [0,\infty) \to \ol{\bD}$ is obtained by 
\begin{align*}
	\eta_{t} = \lim_{r\uparrow 1}\SLEgAt{t}^{-1}(rX_{t}),\quad t \geq 0.
\end{align*}
For the law of the curve to be conformally invariant and satisfy the domain Markov property, the driving process must be $X_{t}=e^{\bbmi \sqrt{\kappa}B_{t}}$, $t\geq 0$
with $\kappa\geq 0$ and a standard Brownian motion $(B_{t}:t\geq 0)$.

It is natural to expect that an analogous story as in the previous section
exists for the radial setting as well.
In fact, the SLE/CFT-correspondence (or coupling) in the radial setting has been studied in~\cite{BauerBernard2004b}. 
As its probabilistic instance,~\cite{IzyurovKytola2013} studied the local set coupling between the radial SLE at $\kappa=4$ and GFF.

\subsection{Summary of the results}
In this paper, we explore the local set coupling between radial multiple SLE and GFF,
thereby extending the results of both~\cite{KK_three_phases_2021} in the chordal setting and~\cite{IzyurovKytola2013} in the single-curve setting.

\begin{figure}[htbp]
    \begin{tabular}{ccc}
%    \hskip -2.5cm
       \centering
        \includegraphics[keepaspectratio, scale=0.3]{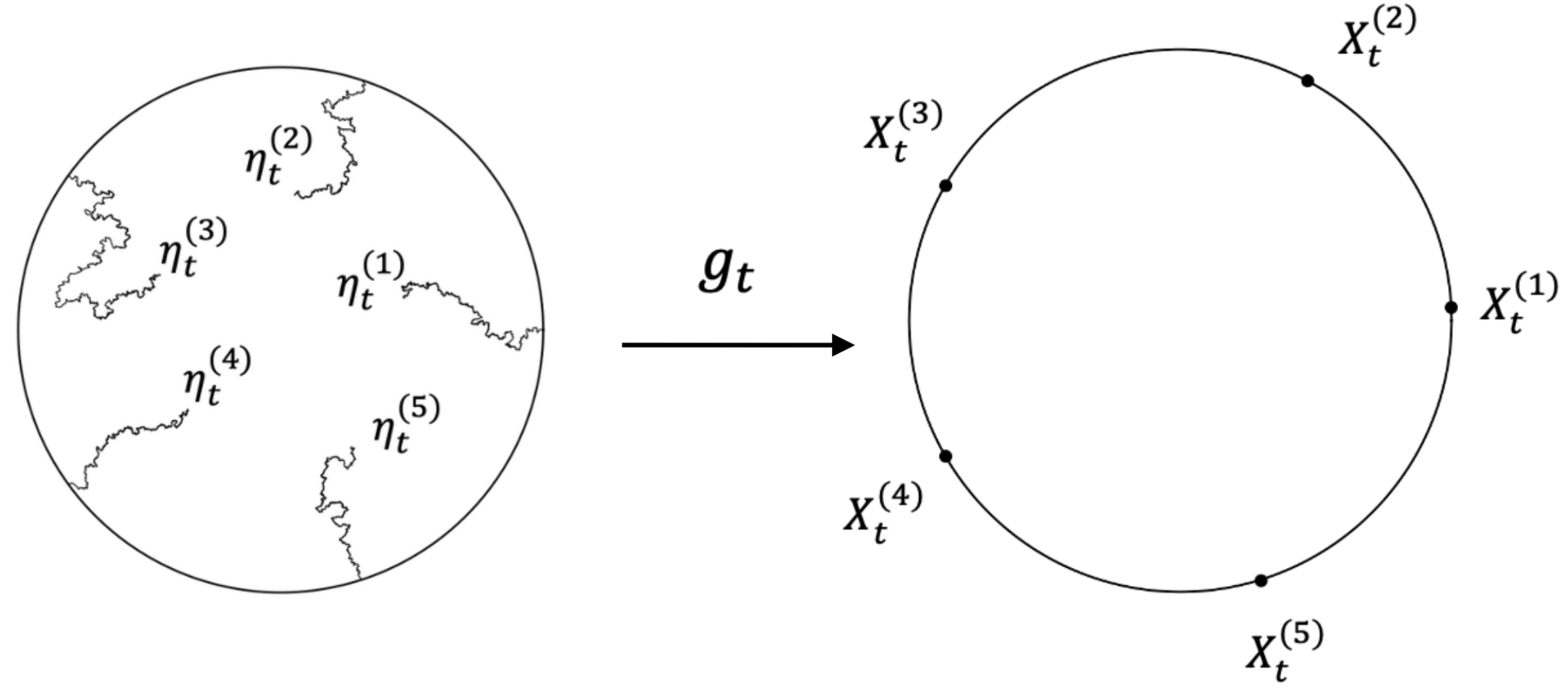}
    \end{tabular}
     \caption{
A schematic picture of 
$\SLEgAt{t} \colon \SLEDomAt{t} \to \bD$ at time $t >0$.
The tips of SLE curves,
$\eta_t^{(i)}$, $i=1,\dots, N$ are
mapped to the points 
$X_t^{(i)}$, $i=1, \dots, N$ on $\bS$,
respectively.
     }
     \label{fig:rmSLE}
  \end{figure}

Our local model for radial multiple SLE 
is the following natural extension of (\ref{eq:radial_loewner}):
with $N\geq 1$ ,
\begin{align}
\label{eq:radial_multiple_loewner}
    \frac{d}{dt}\SLEgAt{t}(z)&=-\sum_{i=1}^{N}\SLEgAt{t}(z)\frac{\SLEgAt{t}+X^{(i)}_{t}}{\SLEgAt{t}-X^{(i)}_{t}},\quad t\geq 0,\\
    \SLEgAt{0}(z)&=z\in \bD, \notag
\end{align}
driven by continuous driving processes $(X^{(i)}_{t}:t\geq 0)$, $i=1,\dots, N$ on $\bS$.
The solution is a family of conformal maps $\SLEgAt{t}\colon \SLEDomAt{t} \to \bD$, $t\geq 0$,
and $N$ curves $\eta^{(i)}\colon [0,\infty) \to \ol{\bD}$, $i=1,\dots, N$ are defined by 
\begin{align*}
\eta^{(i)}_{t}:=\lim_{r\uparrow 1}\SLEgAt{t}^{-1}(rX^{(i)}_{t}), \quad t\geq 0,\quad i=1,\dots, N.
\end{align*}
See Figure~\ref{fig:rmSLE} for illustration.
Note that the curves $\eta^{(i)}$, $i=1,\dots, N$ might not be continuous because the driving processes are yet unknown.
Nevertheless, the domains $\SLEDomAt{t}$, $t\geq 0$ on which the conformal maps $\SLEgAt{t}$, $t\geq 0$ are defined make sense.
For completeness, we write $\SLEHullAt{t}\coloneqq \bD\setminus \SLEDomAt{t}$, $t\geq 0$
and call it the SLE hull.

%%%%%%%%%%%%%%%%%%%%%%%%%%%%%%
%%%%%%%%%%%%%%%%%%%%%%%%%%%%%%%

\begin{defn}
We call the solution of (\ref{eq:radial_multiple_loewner}) the \textit{radial multiple SLE} driven by $(X^{(i)}_{t}:t\geq 0)$, $i=1,\dots, N$.
\end{defn}

Similarly to the chordal case, the law of the driving processes is determined by a partition function~\cite{krusell2024commutation,makarov2025multiple}.
Otherwise, a coupling to GFF provides enough information to fix the driving processes.
To state the result, let us first recall the 
\textit{circular Dyson Brownian motions}.

\begin{defn}
Let $N\geq 1$ and $\beta>0$.
The $N$-particle 
\textit{circular $\beta$-Dyson Brownian motions} 
are $(X^{(i)}_{t}=e^{\bbmi \Theta^{(i)}_{t}}\in \bS: t \geq 0)$, $i=1,\dots,N$ that solve the system of stochastic differential equations (SDEs)
\begin{equation}
\label{eq:circular_Dyson}
    d\Theta^{(i)}_{t}=\sqrt{\frac{8}{\beta}}\,dB^{(i)}_{t}+2\sum_{\substack{j=1\\j\neq i}}^{N}\cot\left(\frac{\Theta^{(i)}_{t}-\Theta^{(j)}_{t}}{2}\right)dt,\quad i=1,\dots, N,
\end{equation}
where $(B^{(i)}_{t}:t\geq 0)$, $i=1,\dots, N$ are mutually independent standard Brownian motions.
\end{defn}

\begin{rem}
The system of SDEs~(\ref{eq:circular_Dyson}) seems non-standard.
In fact, the processes $(X^{(i)}_{\beta t/8}:t\geq 0)$, $i=1,\dots, N$ after changing the time scale are often referred to as the circular $\beta$-Dyson Brownian motions.
The reason for the peculiar coefficient $\sqrt{8/\beta}$ is that we will match $8/\beta$ with the SLE parameter $\kappa$.
\end{rem}

We will phrase out the definition of the local set coupling between radial multiple SLE and GFF in Section~\ref{sect:coupling}.
In our context, a GFF refers to the sum of the zero-boundary GFF and a deterministic harmonic function.
Thus, we need to pick a harmonic function to specify a GFF.
We will find a family of GFFs parametrized by $\kappa>0$ so that the following holds. (See Theorem~\ref{thm:coupling_multi_curve} for precise statement.)

\begin{thm}
\label{thm:coupling_sketch} 
The GFF with parameter $\kappa>0$ is coupled with a radial multiple SLE
if and only if the radial multiple SLE is driven by the circular $8/\kappa$-Dyson Brownian motions.
\end{thm}

After completing the first version of this manuscript, we learned that the ``if'' part of the theorem was independently obtained in~\cite{BerestyckiLisLiuPeltola2025} as well.

%%%%%%%%%%%%%%%%%%%%%%%%%%%%%%
\begin{figure}[htbp]
    \begin{tabular}{ccc}
    \hskip -2.5cm
      \begin{minipage}[t]{0.63\hsize}
       \centering
        \includegraphics[keepaspectratio, scale=0.22]{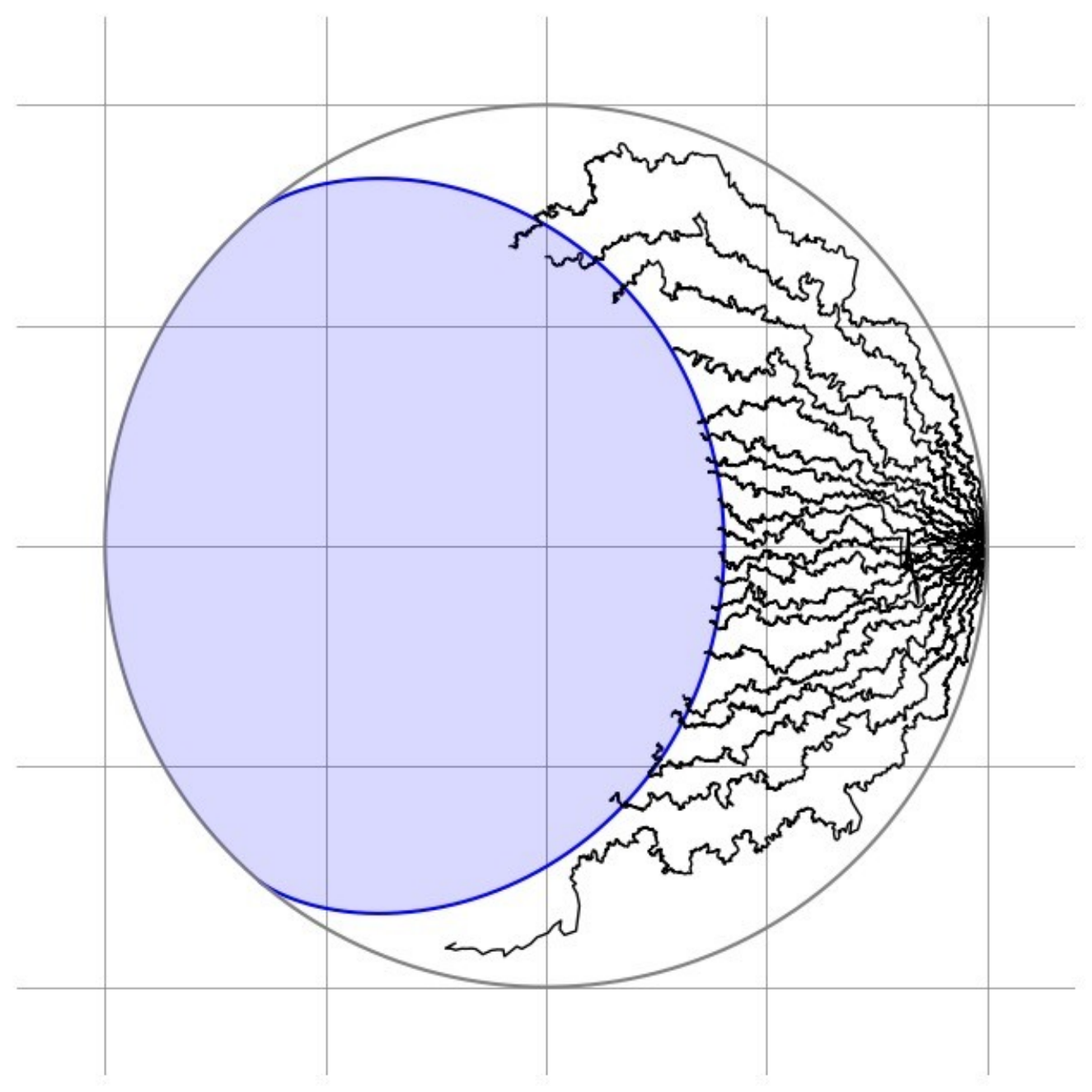}
        \subcaption{}
        \label{fig:SLE_t=0.5}
      \end{minipage} &
    \hskip -5.3cm
      \begin{minipage}[t]{0.63\hsize}
       \centering
        \includegraphics[keepaspectratio, scale=0.22]{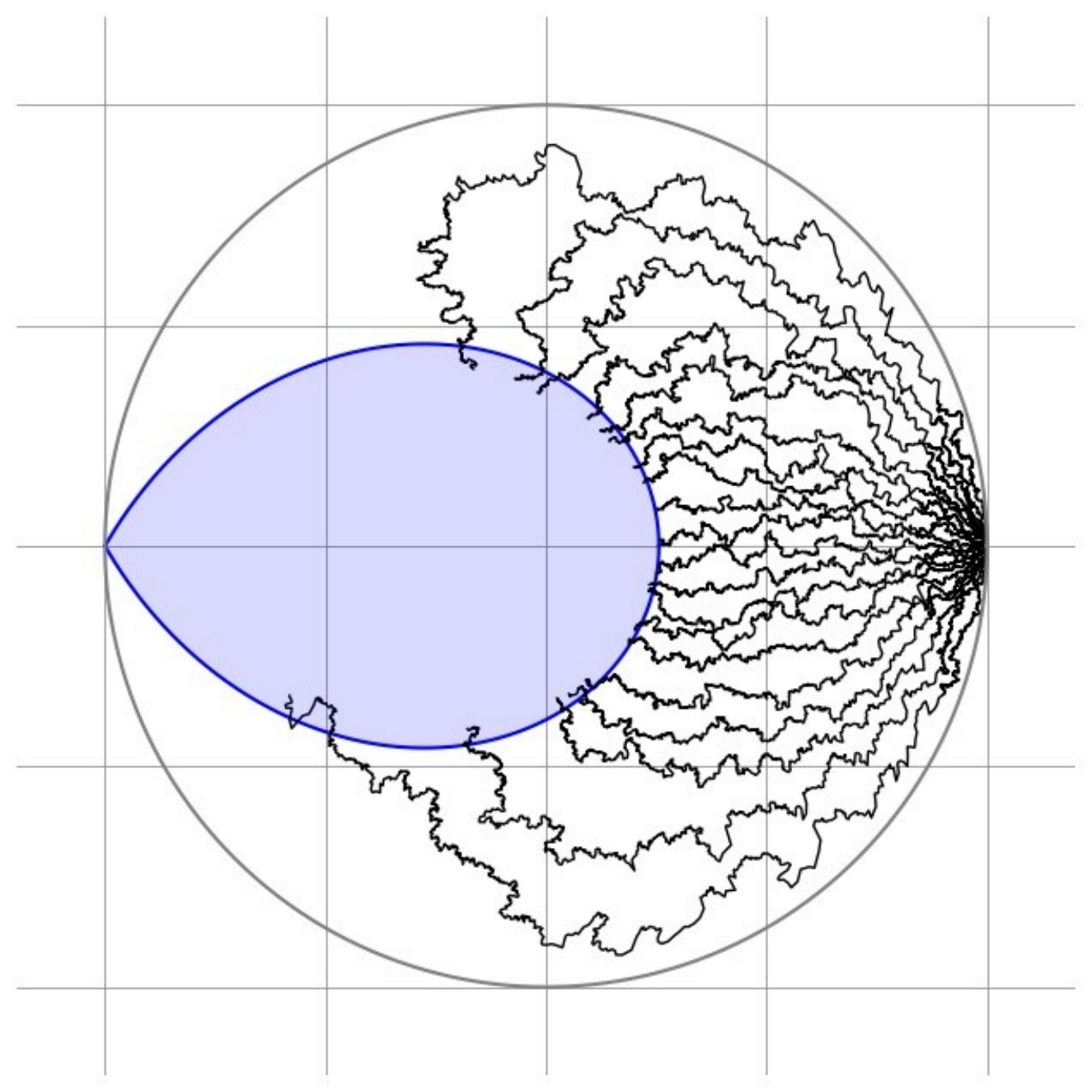}
        \subcaption{}
        \label{fig:SLE_t=1}
      \end{minipage} &
          \hskip -5.3cm
      \begin{minipage}[t]{0.63\hsize}
       \centering
        \includegraphics[keepaspectratio, scale=0.22]{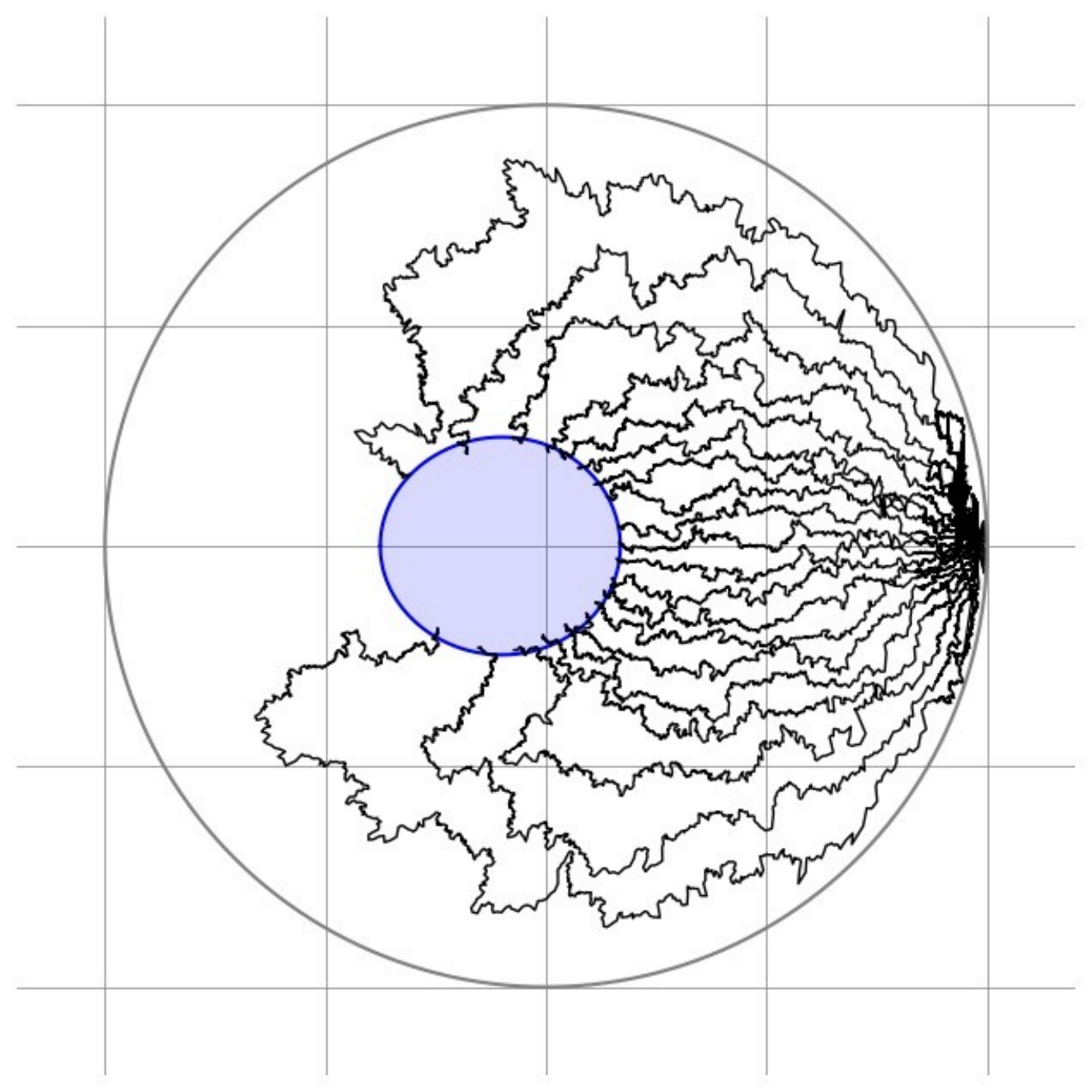}
        \subcaption{}
        \label{fig:SLE_t=1.5}
      \end{minipage} 
    \end{tabular}
     \caption{
Samples of radial SLE curves $\eta_t^{(i)}$, $i=1,\dots, N$
with $N=20$, 
all started from $z=1$, are drawn by numerical 
simulation for $\kappa=8/\beta=2$ at time
(a) $t=0.5$, 
(b) $t=t_{\mathrm{c}}=1$ (critical time), 
(c) $t=1.5$, respectively.
At each time, as $N \to \infty$, 
SLE curves become denser and 
formulate a subdomain of $\D$
which we call the SLE hull $\mSLEHullAt{\infty}{t}$.
The domains of $\mSLEgAt{\infty}{t}$ given by
$\mSLEDomAt{\infty}{t} =\bD \setminus \mSLEHullAt{\infty}{t}$ 
are shown by the shaded subdomains of $\bD$ 
and the boundary curves 
of the SLE hulls $\gamma_t$ are drawn 
by thick lines for each time.
The SLE hull $\mSLEHullAt{\infty}{t}$ changes 
its topology at $t_{\mathrm{c}}$ from a disk to an annulus
}
\label{fig:SLE}
  \end{figure}
%%%%%%%%%%%%%%%%%%%%%%%%%%%%%%%

The circular Dyson Brownian motions are a natural analogue of the Dyson Brownian motions on the unit circle.
The hydrodynamic limit of the circular Dyson Brownian motions is well-known~\cite{cepa2001brownian}.
To explain the hydrodynamic limit of the corresponding radial multiple SLE,
let us first manifest $N$ in our notation.
For each $N\geq 1$, we will write $(\mSLEgAt{N}{t}:t\geq 0)$ 
for the multiple radial SLE driven by the $N$-particle circular $\beta$-Dyson Brownian motions $(X^{(i)}_{t}:t\geq 0)$, $i=1,\dots, N$, with any $\beta>0$, and $(\mSLEDomAt{N}{t}:t\geq 0)$ for the corresponding evolution of domain.

Along with the hydrodynamic limit $N\to \infty$, we assume that the initial conditions for the circular Dyson Brownian motions satisfy
\begin{equation}
\frac{1}{N} \sum_{i=1}^N \diracMeasAt{X^{(i)}_0}
\, \weakConv \, \diracMeasAt{1}
\quad \mbox{as $N \to \infty$}. 
\label{eq:single_source}
\end{equation}
Here, the Dirac (point-mass) measure on $\bS$ is defined by
\begin{align*}
\diracMeasAtOf{x}{A} = 
	\begin{cases}
		1, & x\in A, \\
		0, & \text{otherwise},
	\end{cases}
\end{align*}
for $x\in \bS$ and Borel measurable $A\subset \bS$.

It has been show in~\cite{hotta2021limits} that, 
for each $t\geq 0$, 
the conformal map $\mSLEgAt{N}{t/N}$ 
converges locally uniformly in distribution as $N\to \infty$.
This is equivalent to that 
$\mSLEDomAt{N}{t/N}$ converges in the Caratheodory kernel sense.
Let us write these limits as $\mSLEgAt{\infty}{t}=\lim_{N\to\infty} \mSLEgAt{N}{t/N}$ and $\mSLEDomAt{\infty}{t} = \lim_{N\to\infty}\mSLEDomAt{N}{t/N}$, $t \geq 0$,
and call them the \textit{hydrodynamic limits of
the radial multiple SLE}. 
We also set $\mSLEHullAt{\infty}{t}\coloneqq \bD \setminus \mSLEDomAt{\infty}{t}$, $t\geq 0$,
which is called the \textit{SLE hull in the hydrodynamic limit}. 
Figure~\ref{fig:SLE} shows numerical simulation of 
the radial SLE curves $\eta_t^{(i)}$, $i=1,\dots, N$ 
with $N=20$ for $\kappa=8/\beta=2$ at $t=0.5, 1.0$, and 1.5,
under the initial condition $X_0^{(i)}=\eta_0^{(i)}=1$ for all $i=1,\dots, N$. 
As $N$ grows, the SLE curves will become denser, and
the SLE hull will be eventually a subdomain 
$\mSLEHullAt{\infty}{t}$ of $\D$
such that $e^{\pm \bbmi \varepsilon} \in \ol{\mSLEHullAt{\infty}{t}}$, $0 \leq \varepsilon \ll 1$,
but $0 \notin \mSLEHullAt{\infty}{t}$, $0 < t < \infty$.
The limit family of conformal maps $(\mSLEgAt{\infty}{t}:t\geq 0)$ solves a measure-driven Loewner equation~\cite{hotta2021limits} (see Proposition~\ref{thm:hydro_SLE}) that allows for further analysis.
It is worth noting that the hydrodynamic limit does not depend on the parameter $\beta$ for the circular Dyson Brownian motions.

Here we report exact and explicit description of 
$\mSLEgAt{\infty}{t}$ and $\mSLEDomAt{\infty}{t}$ 
for $t\geq 0$.
Let $W_0$ be the complex Lambert function used in~\cite{HottaKatori2018} (see Section~\ref{sec:Lambert} for definition), and put
\begin{equation}
\Lambda_t(z) =4t \frac{z}{(1-z)^2}.
\label{eq:Lambda}
\end{equation}
Then, the Loewner chain $(\mSLEgAt{\infty}{t}:t\geq 0)$ admits the following expression.
%%%%%%%%%%%%%%%%%%%%%%%%%%
\begin{thm}
\label{thm:hydro_sketch}
At each $t\geq 0$, we have 
\begin{equation*}
\mSLEgAt{\infty}{t}(z) = \left[
1+\frac{2t}{W_0(\Lambda_t(z)e^{-t})}
\left( 1 - \sqrt{1+W_0(\Lambda_t(z)e^{-t})/t} \right)
\right]
e^{2t \sqrt{1+W_0(\Lambda_t(z)e^{-t})/t}}.
\end{equation*}
\end{thm}
%%%%%%%%%%%%%%%%%%%%%%%%%%

To present the result for $\mSLEDomAt{\infty}{t}$, $t\geq 0$, we need a couple of auxiliary functions.
First, we define $x^{\max}_{t}\geq 0$, $t\in [0,\infty)$ and $x^{\min}_{t}\geq 0$, $t\in [1,\infty)$ by the equations
\begin{equation}
\frac{x_t^{\max}}{2}
\tanh \frac{x_t^{\max}}{2}=t \quad \text{and} \quad
\frac{x_t^{\min}}{2}
\coth \frac{x_t^{\min}}{2}=t.
\label{eq:x_min_max}
\end{equation}
It can be readily checked that both are monotonically increasing in $t$,
$x_0^{\max}=x_1^{\min}=0$, 
$x_t^{\min} < x_t^{\max}$ for $t \in [1, \infty)$, and  
$\lim_{t \to \infty} x_t^{\max}=
\lim_{t \to \infty} x_t^{\min}=+\infty$.
Next, we introduce
\begin{equation}
\Phi_t(x) :=
2 \, \arcsinh \left[
- \sqrt{\frac{te^{-t}}{u(x)}}
\exp \left( -\bbmi \psi_{t}(x)-
\frac{u(x)}{2} e^{2 \bbmi \psi_t(x)}
\right) \right]
\label{eq:Phi}
\end{equation}
with
\begin{align}
u(x) &=\frac{x}{\sinh x}, 
\label{eq:ux}
\\
\label{eq:psi}
\psi_t(x) &=\frac{1}{2} \arccos 
\left( \frac{x}{2t} \sinh x-\cosh x \right).
\end{align}
Then, 
\begin{equation}
\label{eq:gamma_t}
\gamma^{+}_t 
= 
\begin{cases}
\big\{ \gamma^{+}_t(x) = e^{\Phi_t(x)} \, \big| \, 
x \in [0, x_t^{\max}] 
\big\},
& \text{if } 0 < t \leq 1,
\\
\big\{ \gamma^{+}_t(x) = e^{\Phi_t(x)} \, \big| \, 
x \in [x_t^{\min}, x_t^{\max}] 
\big\},
& \text{if } t > 1
\end{cases}
\end{equation}
gives a simple curve in $\ol{\bD\cap \bH}$ with the endpoints
\begin{align*}
\gamma_t^{+\mathrm{out}} &:=
\begin{cases}
\gamma^{+}_t(0)=
\begin{cases}
e^{\bbmi \varphi_t^{\mathrm{c}}} \in \bS \cap \bH,
& \text{for }0 < t < 1, \\
-1, & \text{for }t=1,
\end{cases}
\\
\gamma^{+}_t(x_t^{\min}) \in (-1, 0),
\qquad \, \, \, \, \text{for }1< t < \infty,
\end{cases}
\\
\gamma_t^{+\mathrm{in}} &:=\gamma^{+}_t(x_t^{\max}) \in (0, 1),
\qquad \text{for }0 < t < \infty, 
\end{align*}
where
\begin{equation}
\varphi_t^{\mathrm{c}}=
2 \, \arcsin (\sqrt{t} e^{(1-t)/2})
=\arccos (1-2t e^{1-t}).
\label{eq:varphi_t_c}
\end{equation}
Note that
\[
\gamma_t^{+} \setminus
\{\gamma_t^{+\mathrm{out}}, \gamma_t^{+\mathrm{in}} \}
\subset \D \cap \bH.
\]

\begin{rem}
To find out these endpoints, it is useful to notice that $\psi_{t}(x)$ admits various expressions. In fact, we have
\begin{align*}
	\psi_{t}(x) = \arcsin\left(\sqrt{1-\frac{x}{2t}\tanh\frac{x}{2}}\cosh\frac{x}{2}\right),\quad t > 0,\quad x\leq x_{t}^{\max}
\end{align*}
that gives $\psi_{t}(x_{t}^{\max}) = 0$, $t>0$ and $\psi_{t}(0) = \frac{\pi}{2}$, $0<t\leq 1$.
Another formula reads
\begin{align*}
\psi_{t}(x) =\arccos\left(\sqrt{\frac{2}{2t}\coth\frac{x}{2}-1}\sinh\frac{x}{2}\right),\quad t > 1,\quad x\geq x_{t}^{\min}
\end{align*}
that gives $\psi_{t}(x_{t}^{\min}) = \frac{\pi}{2}$, $t>1$.
\end{rem}
For each $t \geq 0$, let 
$\gamma^{-}_t = (\gamma^{+}_t)^{*}$ be
the complex conjugate of \eqref{eq:gamma_t},
and set
\begin{equation}
\gamma_{t} := \gamma^{+}_t \cup \gamma^{-}_t.
\label{eq:SLE_hull}
\end{equation}
We also understand 
$\gamma_0 := \lim_{t \downarrow 0} \gamma^{+}_t
=\{\lim_{t \downarrow 0} \gamma^{+}_t(0) \}
=\{\lim_{t \downarrow 0} \gamma^{+}_t(x_t^{\max})\}
=\{1\}$. 
We can see that $\gamma_{t}$ for $t>0$ cuts $\bD$ into two connected components. 
With these preliminaries, we can identify the domains $\mSLEDomAt{\infty}{t}$, $t\geq 0$ as follows.

\begin{thm}
\label{thm:hydro_domain}
At each $t\geq 0$, $\mSLEDomAt{\infty}{t}$ is the connected component of $\bD \setminus \gamma_{t}$ that contains $0$.
\end{thm}

From the explicit formula for $(\gamma_{t} : t \geq 0)$, 
we can see the following.
There is a \textit{critical time} $t_{\mathrm{c}}=1$
such that 
\begin{enumerate}
\item when $0<t< t_{\mathrm{c}}$, 
$\gamma_{t} = (\partial \mSLEHullAt{\infty}{t}\cap \bD)\cup \{\gamma^{+\mathrm{out}}_t, (\gamma^{+\mathrm{out}}_t)^{*}\}$, 
\item when $t=t_{\mathrm{c}}$, 
$\gamma_{t_{\mathrm{c}}} 
= (\partial \mSLEHullAt{\infty}{t_{\mathrm{c}}} \cap \bD)\cup \{\gamma^{+\mathrm{out}}_{t_{\mathrm{c}}}=(\gamma^{+\mathrm{out}}_{t_{\mathrm{c}}})^{*}=-1\}$,
\item when $t>t_{\mathrm{c}}$, 
$\gamma_{t}=\partial \mSLEHullAt{\infty}{t}\subset \bD$.
\end{enumerate}
In particular, the SLE hull $\mSLEHullAt{\infty}{t}$ changes its topology at time $t_{\mathrm{c}}$ 
from a disk to an annulus as is also demonstrated in Figure~\ref{fig:SLE}.
We will often call $\gamma_{t}$ with $t\geq 0$ the boundary curve of the SLE hull though it does not necessarily coincide with $\partial \mSLEHullAt{\infty}{t}$.

%%%%%%%%%%%%%%%%%%%%%%%%%%%
\begin{rem}
For $0 < t \ll 1$, \eqref{eq:x_min_max} gives
$(x_t^{\max})^2/4 \sim t$.
We parameterize $x \in [0, 2 \sqrt{t}]$ as
$x=2 \sqrt{t} \cos \vartheta$, $\vartheta \in [0, \pi/2]$. 
Equations \eqref{eq:ux} and \eqref{eq:psi} give
$u(x) \sim 1$ and
$\psi_t(x) \sim {\arccos} (2 \cos^2 \vartheta-1)/2
= \vartheta$, respectively.
Hence, the boundary curve of the SLE hull
\eqref{eq:SLE_hull} with \eqref{eq:gamma_t} is reduced to
\begin{equation}
\gamma_t  =\Big\{ 
1-2 \sqrt{t} 
\exp \Big(- \bbmi \vartheta- \frac{1}{2} 
e^{2 \bbmi \vartheta}\Big) \Big| \vartheta \in [- \pi/2, \pi/2] \Big\}.
\label{eq:short}
\end{equation}
If we map the tangent line
to the unit circle $\bS$ at $z=1$ to $\R$ centered at the origin,
then the short-term behavior \eqref{eq:short} is 
identified with the exact expression
of the boundary curve of SLE hull in the hydrodynamic limit
staring from $\delta_0$ in $\bH$, which was
derived in~\cite{HottaKatori2018} for the chordal multiple SLE.
In other words, the present exact solution of
the boundary curve of the SLE hull 
\eqref{eq:SLE_hull} with \eqref{eq:gamma_t} 
can be regarded as a trigonometric-hyperbolic extension of
\eqref{eq:short}. 
\end{rem}

We would like to stress that, even though the results on the hydrodynamic limit is a radial analogue of~\cite{HottaKatori2018} that studied the chordal case, 
their derivation is not exactly parallel to the previous work due to the lack of scaling symmetry.
In fact, it has become manifest in the radial setting that, at the hydrodynamic limit, the analysis of multiple SLE provides information about circular Dyson Brownian motions even though the latter is not explicitly studied.

\subsection{Discussion}
\subsubsection*{Relation to other work}
The circular Dyson Brownian motions were first proposed as the natural driving processes of radial multiple SLE in~\cite{Cardy2003a}.
More recently,~\cite{healey2021n} also suggested that the radial multiple SLE should be driven by the circular Brownian motions.
The starting point of~\cite{healey2021n} is the fact~\cite{KozdronLawler2007} that the global multiple SLE in the chordal setting is obtained by reweighting independent SLEs using Brownian loop measures.
Their requirement on radial multiple SLE is that it is realized as a certain limit of this construction.

Both~\cite{Cardy2003a} and~\cite{healey2021n} assume that, with radial multiple SLE, each curve is absolutely continuous with respect to a single radial SLE,
which is of course reasonable.
It seems to be the case that the driving processes need to be the circular Dyson Brownian motions when all curves must go to the origin.
In contrast, in the present work, we do not assume such absolute continuity, not even that the radial multiple SLE actually generates continuous curves.
Nonetheless, we can fix the driving processes by assuming coupling with GFF.

\subsubsection*{$\beta$-ensemble on the boundary}
The other result of~\cite{Cardy2003a} is that, in a reasonable situation, the radial multiple SLE evolves from random boundary points obeying the circular $\beta$-ensemble, the invariant measure of the circular $\beta$-Dyson Brownian motions.
Though, in the present work, the radial multiple SLE evolves from fixed boundary points, it is interesting to see if GFF realizes a situation where the curves evolves from random boundary points.

In this direction, we find the recent~\cite{liu2025critical} interesting.
The result therein is that, a certain chordal multiple SLE curves coupled with GFF
hit random boundary points in an interval that obeys the $\beta$-Jacobi ensemble.
If a similar approach is possible in the radial case, it could give an answer to the above question.

\subsubsection*{Origin of (circular) $\beta$-Dyson Brownian motions}
Either on the real axis or the unit circle, the (circular) Dyson Brownian motions have a long history. On the real axis, the Dyson Brownian motions for $\beta=1,2,4$ are eigenvalue processes of dynamical random matrices~\cite{Dyson1962}.
The circular Dyson Brownian motions for $\beta=2$ also appeared as an eigenvalue process in the same paper.
At the level of SDEs, it seems natural to extend the parameter to $\beta > 0$ as it is interpreted as the strength of interaction.
In fact, this generalization has led to numerous achievements in the field of stochastic log-gases~\cite{CepaLepingle1997, cepa2001brownian}.

However, once we step back and ask why this generalization to $\beta > 0$ is a right direction, we find the origins of general $\beta$-Dyson Brownian motions quite subtle.
One natural expectation would be that they also appear as eigenvalue processes of some dynamical random matrices.
Potentially, this sounds promising since the static $\beta$-ensemble can be realized as an eigenvalue ensemble of a tridiagonal random matrix~\cite{dumitriu2002matrix}.
However, there turns out to be a discrepancy between one dynamical model of this tridiagonal random matrix that leads to the Dyson Brownian motions~\cite{holcomb2017tridiagonal} and another that does not~\cite{yabuoku2023eigenvalue}.
Another attempt~\cite{allez2012invariant,allez2013diffusive} was restricted to a narrow range of the parameter $\beta\leq 2$.

This subtlety was recently settled and the general $\beta$-Dyson Brownian motions are certainly an eigenvalue process of a dynamical random matrix~\cite{huang2023motion}.
Otherwise, if the Dyson Brownian motions do not have to be an eigenvalue process, we could interpret the previous work~\cite{KK_three_phases_2021} of the first two authors as providing a natural origin of Dyson Brownian motions in terms of the SLE/GFF-coupling.

As far as the authors are aware, the circular $\beta$-Dyson Brownian motions with general $\beta$ have not been realized as an eigenvalue process of a dynamical random matrix though its static counterpart is an eigenvalue ensemble~\cite{dumitriu2002matrix}.
In this regard, we could say that the present work provides a natural origin of circular Dyson Brownian motions similarly to the chordal case.

\subsubsection*{Three phases of radial multiple SLE}
In the previous work~\cite{KK_three_phases_2021} of the first two authors, we applied the coupling result to show that the multiple SLE exhibits three phases against the parameter $\kappa$.
The three phases are well-known for the single curve SLE and can be transferred to the multiple SLE by absolute continuity of each curve with respect to the single curve SLE.
However, the absolute continuity is only manifestly valid as long as the curves are apart from each other, i.e, when $\kappa \leq 4$.
With our approach, we knew that there was a law of multiple curves determined by the GFF, and that allowed us to conclude the three phases for the larger range $\kappa \leq 8$.

We expect that the similar is true for the radial case, but the details are yet to be explored.

\subsubsection*{Other phase transitions}
In the present work, we have seen that, at the hydrodynamic limit, the SLE hull $\mSLEHullAt{\infty}{t}$ changes its topology at $t=1$.
As we will see, this corresponds to the fact that the support of the circular Dyson Brownian motions at the hydrodynamic limit covers the entire $\bS$.
In some sense, we {\it upgraded} a phenomenon for the circular Dyson Brownian motions to the multiple SLE they drive.

There are a few intriguing phase transitions known for the circular Dyson Brownian motions at $\beta=2$, also known as the non-intersecting Brownian motions.
One example is~\cite{forrester2011non} (see also~\cite{majumdar2014top}) that studied the {\it normalized reunion probability} and found a third-order phase transition.
Their method is remarkable; they matched the normalized reunion probability with the free energy of a Yang--Mills theory and transferred 
the known result of phase transition due to~\cite{douglas1993large}.
Thus, their result is a version of the Gross--Witten--Wadia transition~\cite{gross1980possible,wadia1980n} in random matrix theory.

Another phase transition is found in~\cite{liechty2016nonintersecting} that studied the non-intersecting Brownian motions on a circle conditioned to reunion after time $T$.
They found that the widing number for those particles exhibits a critical phenomenon against $T$.

It is interesting to see if these other phase transitions can be {\it upgraded} to the multiple SLE.

\subsubsection*{Beyond hydrodynamic limit}
Now that we understand the hydrodynamic limit of multiple SLE,
a natural next question would be to see fluctuation in different scaling.
The first result beyond the hydrodynamic limit appears in~\cite{campbell2023rate} that studies the convergence rate in the chordal setting. 
We expect a similar result to hold in the radial case as well.

In random matrix theory, the Tracy--Widom distribution~\cite{tracy1994level} for the properly scaled largest eigenvalue has been playing a prominent role in revealing universality.
It would be interesting if we can lift the Tracy--Widom distribution to the multiple SLE, for example, by analysing the distribution of the right-most curve in a particular scaling.

Even further, we could go on to the large deviations.
In fact, the third-order phase transition found in~\cite{forrester2011non,majumdar2014top} corresponds to the discrepancy between the left and right rate functions away from the critical value.
Therefore, studying large deviations would be relevant for answering the questions posed in the above subsection.
Also, the large deviation principles for the empirical measure of Dyson Brownian motions have been established in~\cite{guionnet2002large,guionnet2021large}.
Although it is unclear what is the analogue of the empirical measure for multiple SLE, it is interesting if multiple SLE admits a certain large deviation principle.

Though orthogonal to our setting of sending $N\to\infty$, we also would like to mention~\cite{peltola2023large,abuzaid2024large} that study the large deviation principles of multiple SLE as $\kappa\to 0$.

\subsubsection*{Multiply connected domains}
The notion of SLE has been generalized to multiply connected domains~\cite{Zhan2004,lawler2006laplacian, bauer2006radial,BauerFriedrich2008,drenning2011excursion,chen2017stochastic,ChenFukushima2018,murayama2019chordal}.
On the other hand, the definition of GFF does not require that the domain is simply connected from the first place.
In fact, coupling between these variants of SLE and GFF on multiply connected domains appeared in~\cite{IzyurovKytola2013} on annuli and has been recently studied more in~\cite{byun2023conformal,alberts2024conformal,aru2024sle}.
The present work encourages us to study the coupling in presence of multiple curves and derive a natural stochastic log-gas that arises in the boundaries of a multiply connected domain.
This direction will be pursued in future work.

\subsection*{Organization of the paper}
This paper is organized as follows.
The following Section~\ref{sect:GFF} provides a brief exposition of GFF.
In Section~\ref{sect:coupling}, we discuss the local set coupling between radial multiple SLE and GFF. This section contains the definition of the coupling, a precise version of Theorem~\ref{thm:coupling_sketch} as well as its proof.
Finally, in Section~\ref{sect:hydro}, we describe the hydrodynamic limit of the radial multiple SLE driven by the circular Dyson Brownian motions.

\subsection*{Acknowledgments}
MK, CS, and RS thank Taiki Endo and Miyu Ikegame for 
collaboration in the early stages of the study 
reported in Section \ref{sect:hydro}.
The authors are grateful to Sebastian Schlei{\ss}inger for informing them of the GitHub repository~\url{https://github.com/Sammy-Jankins/Complex_analysis}, 
from which they adapted Python codes, as well as for comments on the classical limit.
They also thank Saori Morimoto for useful discussion on numerical simulations that produced Figures~\ref{fig:rmSLE} and~\ref{fig:SLE}.
MK was supported by JSPS KAKENHI Grant Numbers 
JP21H04432,
JP22H05105,
JP23K25774,
and
JP24K06888.
SK was supported by Research Council of Finland, project 340965: {\it Coupling of Gaussian free fields and Schramm–Loewner evolutions on multiply connected domains} and project 346309: {\it Finnish Centre of Excellence in Randomness and Structures, ``FiRST''}.

\section{Gaussian free field}
\label{sect:GFF}
In this section, we introduce the necessary background on Gaussian free field (GFF). For more detailed account, see~\cite{Sheffield2007,werner2020lecture, berestycki2024gaussian}.

\subsection{Definition}
Let $D\subset\bC$ be a domain that has the Green's function $G_{D}(z,w)$ under the zero-boundary condition.
Our convention is that the Green's function has the singularity
\begin{equation*}
    G_{D}(z,w)\sim -\log |z-w|\quad \mbox{as}\, z\to w.
\end{equation*}
Under this normalization, $\frac{-1}{2\pi}G_{D}(z,w)$ is the integral kernel of the inverse of the Laplacian $\Delta^{-1}$:
\begin{equation}
\label{eq:inverse_Laplacian_Green_function}
    (-\Delta^{-1}f)(z)=\frac{1}{2\pi}\int_{D}G_{D}(z,w)f(w)\leb(dw),
\end{equation}
where $\leb(dw):=d(\rmRe w) d(\rmIm w)$ is the Lebesgue measure on $D$ induced from $\bC$.

\begin{exam}
\label{exam:Green_D}
We will need the Green's function on $\bD$.
It is given by
\begin{equation*}
    G_{\bD}(z,w)=\log\Bigl|\frac{z\ol{w}-1}{z-w}\Bigr|,\quad z,w\in \bD.
\end{equation*}
Indeed, it has the correct singularity at $z=w$, is symmetric under exchanging $z$ and $w$, and when $w\in \bS$, we can see that $G_{\bD}(z,w)=0$ since $\ol{w}=w^{-1}$.
\end{exam}

We write $C^{\infty}_{0}(D)$ for the space of smooth functions on $D$ with compact support.

\begin{defn}
The zero-boundary GFF on $D$ is a random distribution $\gffD$ with test functions in $C^{\infty}_{0}(D)$ such that
$\gffD (f)$, $f\in C^{\infty}_{0}(D)$ are centered Gaussian variables
with covariance given by
\begin{align*}
	\Cov (\gffD (f), \gffD (g)) = \int_{D\times D}f(z)G_{D}(z,w)g(w)\leb(dz)\leb(dw),\quad f,g\in C^{\infty}_{0}(D).
\end{align*}
\end{defn}

A non-trivial question is if a GFF even exists.
We shall argue that it exists by sketching its construction.
We start off by endowing $C^{\infty}_{0}(D)$ with the Dirichlet inner product
\begin{equation*}
    (f,g)_{\nabla}=\frac{1}{2\pi}\int_{D} (\nabla f)(z) \cdot (\nabla g)(z) \leb(dz),\quad f,g\in C^{\infty}_{0}(D).
\end{equation*}
The Hilbert space completion of $C^{\infty}_{0}(D)$ with respect to this inner product will be denoted by $W(D)$.
For a complete orthonormal system $\{\consAt{i}\}_{i}$ of $W(D)$, we consider the following infinite sum
\begin{equation}
\label{eq:gff_infinite_sum}
    \gffD =\sum_{i}\alpha_{i}\consAt{i},\quad \alpha_{i}\sim N(0,1):\mbox{i.i.d.}
\end{equation}
It is clear that $\gffD$ converges in $W(D)$ with zero-probability.
The hardest part of the construction of a GFF is that $\gffD$ almost surely converges as a distribution with test functions in $C^{\infty}_{0}(D)$.
We do not give a complete proof of this part, but just mension that it relies on Weyl's law for the spectrum of the Laplacian.
See~\cite{Sheffield2007} for detail.

Once we know that $\gffD$ converges as a distribution, it is clear that $\gffD (f)$, $f\in C^{\infty}_{0}(D)$ are centered Gaussian variables.
It remains to see that they exhibit the correct covariance.

\begin{prop}
The covariance of the centered Gaussian variables $\gffD (f)$, $f\in C^{\infty}_{0}(D)$ are given by
\begin{equation*}
    \Cov (\gffD (f),\gffD (g))=\int_{D\times D}f(z)G_{D}(z,w)g(w)\leb(dz)\leb(dw),\quad f,g\in C^{\infty}_{0}(D).
\end{equation*}
\end{prop}
\begin{proof}
Since the assignment $f\mapsto \gffD (f)$ is linear, it suffices to check that each $\gffD (f)$, $f\in C^{\infty}_{0}(D)$ exhibits the variance
\begin{equation*}
    \Var (\gffD (f))=\int_{D\times D}f(z)G_{D}(z,w)f(w)\leb(dz)\leb(dw).
\end{equation*}
Since $\{\consAt{i}\}_{i}$ is orthonormal with respect to the Dirichlet inner product, by integration by parts, we can see
\begin{equation*}
    \delta_{i,j}=(\consAt{i},\consAt{j})_{\nabla}=\frac{-1}{2\pi}(\consAt{i},\Delta\consAt{j}),
\end{equation*}
which implies that $\{(\frac{-1}{2\pi}\Delta)\consAt{i}\}_{i}$ is the dual basis of $\{\consAt{i}\}_{i}$ with respect to the $L^{2}$-pairing.
Let us expand $f\in C^{(\infty)}_{0}(D)$ in terms of $\{(\frac{-1}{2\pi}\Delta)\consAt{i}\}_{i}$:
\begin{equation*}
    f=\sum_{i}f_{i}\left(\frac{-1}{2\pi}\Delta\right)\consAt{i}.
\end{equation*}
Then, we obtain
\begin{equation*}
    \gffD (f)=\sum_{i}\alpha_{i}f_{i},
\end{equation*}
from which we get
\begin{align*}
    \Var (\gffD (f))=\sum_{i}f_{i}^{2}.
\end{align*}
By definition, each $f_{i}$ is given by
\begin{equation*}
    f_{i}=(f,\consAt{i})=((-2\pi\Delta^{-1})f,\consAt{i})_{\nabla},
\end{equation*}
which leads to
\begin{equation*}
    \Var (\gffD (f))=\|(-2\pi\Delta^{-1})f\|_{\nabla}^{2}.
\end{equation*}
By integration by parts, we obtain the desired result:
\begin{equation*}
    \Var (\gffD (f))=2\pi \int_{D} ((-\Delta)^{-1}f)(z)f(z)dz=\int_{D\times D}f(z)G(z,w)f(w)\leb(dz)\leb(dw),
\end{equation*}
where we also used (\ref{eq:inverse_Laplacian_Green_function}).
\end{proof}

More generally, we refer to the sum of the zero-boundary GFF and a deterministic harmonic function as a GFF.
Since a harmonic function is determined by its boundary values, this is to impose boundary conditions on the GFF. 

\subsection{Domain Markov property}
Let us suppose that $D\subset \bC$ is a domain as above and $\gffD_{D}$ is a Dirichlet boundary GFF on $D$. Note that, since $\gffD_{D}$ is a distribution on $D$, we can restrict it to a subdomain.
The domain Markov property of GFF refers to the following property.

\begin{prop}
\label{prop:gff_domain_markov_property}
Let us pick a subdomain $U\subset D$. Then, conditioned over $D\setminus U$,
\begin{align*}
	\gffD_{D}|_{U} = \gffD_{U} + \harm,
\end{align*}
where $\harm$ is the harmonic extension of $\gffD_{D} |_{\dee U}$ in $U$.
\end{prop}

We only give heuristics instead of a complete proof.
Note that, when we embed $W(U)$ in $W(D)$, the orthogonal complement $W(U)^{\perp}$ consists of those in $W(D)$ that are harmonic on $U$.
We may pick complete orthonormal systems of $W(U)$ and $W(U)^{\perp}$
and combine them to get one of $W(D)$.
Therefore, the infinite-sum (\ref{eq:gff_infinite_sum}) decomposes into
mutually independent parts as
\begin{align*}
	\gffD_{D} = \gffD_{U} + \gffD^{h}_{U},
\end{align*}
where $\gffD_{U}$ is a Dirichlet boundary GFF on $U$, extended as zero outside $U$, and $\gffD^{h}_{U}$ is harmonic on $U$.
When conditioned over $D\setminus U$ and restricted on $U$, $\gffD_{U}$ remains as is and $\harm = \gffD^{h}_{U}|_{U}$ gives a harmonic extension of $\gffD_{D}|_{\dee U}$ on $U$.

From the above explanation, we can see that the harmonic function $\harm$ in Proposition~\ref{prop:gff_domain_markov_property} can be extracted as a conditional expectation.
We write $\cF_{(D\setminus U)^{+}}$ for the $\sigma$-algebra generated by $\gffD$ projected to $W(U)^{\perp}$.
Since $\gffD_{U}$ is centered, we have 
\begin{align*}
	\harm = \bE[\gffD_{D}|_{U}|\cF_{(D\setminus U)^{+}}].
\end{align*}
This observation will be used later.

\section{Coupling}
\label{sect:coupling}
In this section, we discuss the local set coupling between radial multiple SLE and GFF,
and present our first main result, Theorem~\ref{thm:coupling_multi_curve},
the precise version of Theorem~\ref{thm:coupling_sketch}.

\subsection{Motivating example: single-curve case}
The local set coupling between radial SLE and GFF has been studied in~\cite{IzyurovKytola2013} in the case that the SLE parameter is $\kappa=4$.
To motivate our setting with radial multiple SLE, let us first formulate the result for the single-curve case.
(This would also be useful as computational warm-up.)

Let $\kappa>0$ and $(g_{t}:t\geq 0)$ be the radial SLE driven by $(X_{t}=e^{\bbmi \sqrt{\kappa} B_{t}}:t\geq 0)$, where $(B_{t}:t\geq 0)$ is a standard Brownian motion, i.e., it is the solution to (\ref{eq:radial_loewner}).
Recall that, at each $t\geq 0$, we may find a domain $\bD_{t}\subset \bD$ that is conformally mapped to $\bD$ by $g_{t}$.
At each $t\geq 0$, we define the harmonic function $\harmAt{t}$ on $\bD_{t}$ by
\begin{align}
\label{eq:harm_fun_single_curve}
    \harmAt{t}(z)=-\frac{2}{\sqrt{\kappa}}\arg (g_{t}(z)-X_{t})+\xi_{1}\arg (g_{t}(z))-\chi \arg (g'_{t}(z))+B_{t},\quad z\in \bD_{t},
\end{align}
where $\xi_{1}=\frac{3}{\sqrt{\kappa}}-\frac{\sqrt{\kappa}}{2}$ and $\chi=\frac{2}{\sqrt{\kappa}}-\frac{\sqrt{\kappa}}{2}$.
Note that $\harmAt{t}$ is not single valued due to the presence of $\arg (g_{t}(z))$ and $\arg (g'_{t}(z))$.
Nevertheless, we argue that $\harmAt{t}$ is the natural generalization of the harmonic function that was used to construct a coupling of chordal SLE and GFF.

Let $U\subset \bD$ be an open set and define the stopping time
\begin{align}
\label{eq:stopping_time_U}
	\tau_{U} = \sup\{t\geq 0|U\subset \bD_{t}\}.
\end{align}
For a stopping time $\tau$ such that $\bP[\tau\leq \tau_{U}]=1$, we define the GFF $\gffB_{U,\tau}$ on $U$ by first sampling $\bD_{\tau}$ and then setting
\begin{align}
\label{eq:gff_single_curve_cut}
	\gffB_{U,\tau} = (\gffD_{\bD_{\tau}}+\harm_{\tau})|_{U}.
\end{align}

In informal terms, the local set coupling between SLE and GFF refers to the invariance of GFF under the time evolution according to SLE.
The following theorem has been proved in~\cite{IzyurovKytola2013} in the case that $\kappa=4$.
\begin{thm}
\label{thm:coupling_single_curve}
For any open set $U\subset \bD$ and a stopping time $\tau$ such that $\bP[\tau\leq \tau_{U}]=1$, we have 
\begin{align*}
	\gffB_{U,\tau} \overset{\law}{=} \gffB_{U,0}.
\end{align*}
\end{thm}

\begin{rem}
The last term $B_{t}$ in (\ref{eq:harm_fun_single_curve}) is missing from the harmonic function used in~\cite{IzyurovKytola2013}.
This is necessary for the coupling in Theorem~\ref{thm:coupling_single_curve} to our understanding.
\end{rem}

\begin{rem}
The terminology came from the fact that $K_{\tau}=\bD\setminus\bD_{\tau}$ is a local set for $\gffD_{\bD}+\harm_{0}$ (see~\cite{MillerSheffield2016a} for the definition).
Although the literal local set coupling is stronger than Theorem~\ref{thm:coupling_single_curve}, we shall not fill this gap as the proof would be very similar to the chordal case~\cite{MillerSheffield2016a}, provided Theorem~\ref{thm:coupling_single_curve}.
\end{rem}

\begin{proof}[Proof of Theorem~\ref{thm:coupling_single_curve}]
First, we show that, for each $z\in \bD$, the process $(\harmAt{t}(z):0\leq t\leq \tau_{z})$ is a local martingale.
Here, $\tau_{z}$ is the stopping time until which the solution $g_{t}(z)$ exists.
For computational convenience, we shall note that $\harmAt{t}(z)=\rmIm \wtilde{\harm}_{t}(z)$ with
\begin{equation*}
    \wtilde{\harm}_{t}(z)=-\frac{2}{\sqrt{\kappa}}\log (g_{t}(z)-X_{t})+\xi_{1}\log (g_{t}(z))-\chi \log (g'_{t}(z)) + \bbmi B_{t}.
\end{equation*}

Let us execute Ito calculus. The circular Brownian motion $(X_{t}:t\geq 0)$ satisfies
\begin{equation*}
    dX_{t}=\bbmi \sqrt{\kappa}X_{t}dB_{t}-\frac{\kappa}{2}X_{t}dt.
\end{equation*}
Introducing the function
\begin{equation}
\label{eq:radial_Loewner_vect_field}
    \Psi (z,x)=-z\frac{z+x}{z-x},\quad z\in \bD,\, x\in \bS,
\end{equation}
the radial Loewner equation reads
\begin{equation*}
    \frac{d}{dt}g_{t}(z)=\Psi (g_{t}(z),X_{t}).
\end{equation*}
It will be computationally useful to have the Laurent expansion of $\Psi (z,x)$ at $z=x$:
\begin{equation}
\label{eq:radial_Loewner_vect_field_Laurent}
    \Psi (z,x)=-\frac{2x^{2}}{z-x}-3x-(z-x).
\end{equation}
Let us start by computing the increment of $\log (g_{t}(z)-X_{t})$:
\begin{align*}
    d\log (g_{t}(z)-X_{t})=&\frac{d(g_{t}(z)-X_{t})}{g_{t}(z)-X_{t}}-\frac{1}{2}\frac{d\braket{X_{t},X_{t}}}{(g_{t}(z)-X_{t})^{2}} \\
    =&\frac{-\bbmi \sqrt{\kappa}X_{t}}{g_{t}(z)-X_{t}}dB_{t}+K_{1}(g_{t}(z),X_{t})dt,
\end{align*}
where
\begin{align*}
    K_{1}(z,x)=&\frac{\Psi (z,x)}{z-x}+\frac{\kappa}{2}\frac{x}{z-x}+\frac{\kappa}{2}\frac{x^{2}}{(z-x)^{2}} \\
    =&\frac{(\frac{\kappa}{2}-2)x^{2}}{(z-x)^{2}}+\frac{(\frac{\kappa}{2}-3)x}{z-x}-1.
\end{align*}
The second and third terms of $\wtilde{\harm}_{t}(z)$ are purely drift terms: we obtain
\begin{align*}
    \frac{d}{dt}\log (g_{t}(z))=\frac{1}{g_{t}(z)}\frac{d}{dt}g_{t}(z)=L_{1}(g_{t}(z),X_{t})
\end{align*}
with
\begin{equation*}
    L_{1}(z,x)=\frac{1}{z}\Psi (z,x)=-\frac{z+x}{z-x}=-\frac{2x}{z-x}-1
\end{equation*}
and
\begin{align*}
    \frac{d}{dt}\log (g'_{t}(z))=\frac{1}{g'_{t}(z)}\frac{d}{dt}g'_{t}(z)=\frac{\dee \Psi}{\dee z}(g_{t}(z),X_{t})=M_{1}(g_{t}(z),X_{t})
\end{align*}
with
\begin{equation*}
    M_{1}(z,x)=\frac{\dee \Psi}{\dee z}(z,x)=\frac{2x^{2}}{(z-x)^{2}}-1.
\end{equation*}

We shall observe that
\begin{align*}
    & -\frac{2}{\sqrt{\kappa}}K_{1}(z,x)+\xi_{1}L_{1}(z,x)-\chi M_{1}(z,x) \\
    =&\frac{-2\bigl(\chi-(\frac{2}{\sqrt{\kappa}}-\frac{\sqrt{\kappa}}{2})\bigr)x^{2}}{(z-x)^{2}}+\frac{-2\bigl(\xi_{1}-(\frac{3}{\sqrt{\kappa}}-\frac{\sqrt{\kappa}}{2})\bigr)x}{z-x}+\frac{2}{\sqrt{\kappa}}-\xi_{1}+\chi
\end{align*}
is a real constant if $\chi=\frac{2}{\sqrt{\kappa}}-\frac{\sqrt{\kappa}}{2}$ and $\xi_{1}=\frac{3}{\sqrt{\kappa}}-\frac{\sqrt{\kappa}}{2}$.
Therefore, under these choices of $\xi_{1}$ and $\chi$, $(\harmAt{t}(z):0\leq t\leq \tau_{z})$ is a local martingale:
\begin{align*}
    d\harmAt{t}(z)=&-\frac{2}{\sqrt{\kappa}}\rmIm \left(\frac{-\bbmi \sqrt{\kappa}X_{t}}{g_{t}(z)-X_{t}}\right)dB_{t}+dB_{t}=\rmIm\left(\bbmi \frac{g_{t}(z)+X_{t}}{g_{t}(z)-X_{t}}\right)dB_{t} \\
   =&\rmRe \left(\frac{g_{t}(z)+X_{t}}{g_{t}(z)-X_{t}}\right)dB_{t}.
\end{align*}

Next, we pick $z,w\in \bD$ and relate the quadratic variation $\braket{\harmAt{t}(z),\harmAt{t}(w)}$ to the Green's function of $\bD_{t}$.
From the above computation, we have 
\begin{align*}
	\frac{d}{dt}\braket{\harmAt{t}(z),\harmAt{t}(w)} = \rmRe \left(\frac{g_{t}(z)+X_{t}}{g_{t}(z)-X_{t}}\right)\rmRe \left(\frac{g_{t}(w)+X_{t}}{g_{t}(w)-X_{t}}\right),\quad 0\leq t\leq \tau_{z} \wedge \tau_{w}.
\end{align*}
For the Green's function, we have $G_{\bD_{t}}(z,w) = G_{\bD}(g_{t}(z),g_{t}(w))$, and from Example~\ref{exam:Green_D},\begin{align*}
    \frac{d}{dt}G_{\bD_{t}}(z,w)=\Pi_{1}(g_{t}(z),g_{t}(w);X_{t})
\end{align*}
with
\begin{align}
\label{eq:Green_Loewner_variation}
    \Pi_{1}(z,w;x)=&\rmRe \frac{-1}{1-z\ol{w}}\left(\Psi(z,x)\ol{w}+z\ol{\Psi(w,x)}\right)+\rmRe \frac{-1}{z-w}\left(\Psi(z,x)-\Psi(w,x)\right).
\end{align}
By direct computation, we can see that
\begin{equation*}
    \Psi (z,x)\ol{w}+z\ol{\Psi (w,x)}=(1-z\ol{w})\frac{2z\ol{w}}{(z-x)(\ol{w}-\ol{x})},
\end{equation*}
where we used $|x|=1$ for $x\in \bS$ and
\begin{equation*}
    \Psi(z,x)-\Psi(w,x)=(z-w)\left(\frac{z+x}{z-x}\frac{w+x}{w-x}-\frac{2zw}{(z-x)(w-x)}\right).
\end{equation*}
We can immediately notice that, in $\Pi_{1} (z,w;x)$, there is a contribution of the form
\begin{align*}
    \rmRe \left( \frac{z}{z-x}\left(\frac{2w}{w-x}-\frac{2\ol{w}}{\ol{w}-\ol{x}}\right)\right)
    =&\, \rmRe \left(\bbmi \frac{2z}{z-x}\right) \rmIm \left(\frac{2w}{w-x}\right) \\
    =&-\rmIm\left(\frac{2z}{z-x}\right) \rmIm \left(\frac{2w}{w-x}\right).
\end{align*}
Hence, we obtain
\begin{align*}
    \Pi_{1} (z,w;x)=-\rmIm\left(\frac{2z}{z-x}\right)\rmIm \left(\frac{2w}{w-x}\right)-\rmRe\left(\frac{z+x}{z-x}\frac{w+x}{w-x}\right).
\end{align*}
As far as the imaginary part is concerned, we have
\begin{equation*}
    \rmIm \left(\frac{2z}{z-x}\right)=\rmIm \left(\frac{z+x}{z-x}\right).
\end{equation*}
Thus, we get
\begin{equation}
\label{eq:Hadamard_formula}
    \Pi_{1} (z,w;x)=-\rmRe \left(\frac{z+x}{z-x}\right)\rmRe \left(\frac{w+x}{w-x}\right),
\end{equation}
which allows us to conclude that 
\begin{align}
\label{eq:cross_var_harm}
	\braket{\harmAt{t}(z),\harmAt{t}(w)} + G_{\bD_{t}}(z,w) = G_{\bD}(z,w),\quad 0\leq t\leq \tau_{z}\wedge \tau_{w}. 
\end{align}

Now, we pick an open set $U\subset \bD$ and a test function $f\in C^{\infty}_{0}(\bD)$ such that $\supp(f)\subset U$, we have 
\begin{align*}
	\braket{(\harmAt{t},f),(\harmAt{t},f)} + E_{\bD_{t}}(f) = E_{\bD}(f),\quad 0\leq t\leq \tau_{U},
\end{align*}
where 
\begin{align*}
	E_{\bD_{t}}(f) = \int_{\bD_{t}\times \bD_{t}}f(z)G_{\bD_{t}}(z,w)f(w)\leb(dz)\leb(dw)
\end{align*}
is the Dirichlet energy of $f$ in $\bD_{t}$.
Though this is a direct consequence of (\ref{eq:cross_var_harm}) and Fubini's theorem, applying Fubini's theorem requires a slight care. We do not provide the details here as the argument is the same as~\cite[Section~3.3]{MillerSheffield2016a}.

Let us write $(\cF_{t})_{t\geq 0}$ for the filtration to which the radial SLE is adapted.
For any stopping time $\tau$ such that $\bP[\tau\leq\tau_{U}]=1$, $\harm_{\tau}$ is measurable with respect to $\cF_{\tau}$.
Therefore, we have
\begin{align*}
	\bE \bigl[e^{\bbmi \theta (\gffB_{U,\tau},f)}\bigr] =& \bE \biggl[\bE \bigl[e^{\bbmi \theta(\gffD_{\bD_{\tau}},f)}|\cF_{\tau}\bigr]e^{\bbmi \theta (\harm_{\tau},f)}\biggr] \\
	=&\bE \bigl[e^{\bbmi \theta (\harm_{\tau},f)-\frac{\theta^{2}}{2}E_{\bD_{\tau}}(f)}\bigr] \\
	=&\bE \bigl[e^{\bbmi \theta (\harm_{0},f)-\frac{\theta^{2}}{2}E_{\bD}(f)}\bigr] \\
	=&\bE \bigl[e^{\bbmi \theta (\gffB_{U,0},f)}\bigr],\quad \theta\in \bR,
\end{align*}
as is desired.
\end{proof}

\subsection{Multiple curve case}
We now consider the local set coupling of radial multiple SLE and GFF.
Let $N\geq 1$ and $(g_{t}:t\geq 0)$ be the radial multiple SLE with continuous driving processes $(X^{(i)}_{t}=e^{\bbmi \Theta^{(i)}_{t}}:t\geq 0)$, $i=1,\dots, N$, i.e., it is the solution of (\ref{eq:radial_multiple_loewner}).
In contrast to the previous section, we have not fixed the law of the driving processes.
We write $T$ for the collision time of the driving processes that could be finite with positive probability.

To formulate our result, we first need to postulate an analogue of (\ref{eq:harm_fun_single_curve}) when there are multiple driving processes.
The harmonic function \eqref{eq:harm_fun_single_curve} imposes the boundary condition that has a jump at $\SLEgAt{t}(z)=X_{t}$ by $-\frac{2}{\sqrt{\kappa}}$ and an additive constant fluctuation by $B_{t}$.
In the presence of multiple driving processes $(X^{(i)}_{t}:t\geq 0)$, $i=1,\dots, N$,
it would be natural if the boundary value jumps at every $\SLEgAt{t}(z)=X^{(i)}_{t}$, $i=1,\dots, N$.
For the constant fluctuation, we postulate $\sum_{i=1}^{N}\Theta^{(i)}_{t}$.
That is, we fix $\kappa>0$, $\xi,\chi,\zeta\in\bR$, and, at each $t\geq 0$, we take the harmonic function
\begin{align*}
    \harmMAt{N}{t}(z)=& -\frac{2}{\sqrt{\kappa}}\sum_{i=1}^{N}\arg (g_{t}(z)-X^{(i)}_{t})+\xi \arg (g_{t}(z))-\chi \arg (g'_{t}(z))
    +\zeta\sum_{i=1}^{N}\Theta^{(i)}_{t}
\end{align*}
of $z\in \bD_{t}$.
Note that the four parameters $\kappa, \xi, \chi$ and $\zeta$ are not yet related.

Next, we expand the definition of the local set coupling. For an open set $U\subset \bD$, the stopping time $\tau_{U}$ is given by 
\begin{align*}
	\tau_{U} = \sup \{t\geq 0| U\subset \bD_{t}\} \wedge T.
\end{align*}
This is almost same as (\ref{eq:stopping_time_U}) except we cannot continue the process after collision of the driving processes at time $T$.
For a stopping time $\tau$ such that $\bP [\tau \leq \tau_{U}]=1$, the GFF $\gffB_{U,\tau}$ is sampled in the same way as (\ref{eq:gff_single_curve_cut}) except that we use the harmonic function $\harmMAt{N}{\tau}$ instead.
To be precise, we first sample $\bD_{\tau}$, and put 
\begin{align*}
	\gffB_{U,\tau} = (\gffD_{\bD_{\tau}}+\harmMAt{N}{\tau})|_{U}.
\end{align*}

The main result in this section goes as follows.
\begin{thm}
\label{thm:coupling_multi_curve}
The following are equivalent:
\begin{enumerate}
\item 	for any open set $U\subset \bD$ and a stopping time $\tau$ such that $\bP[\tau\leq \tau_{U}]=1$, we have 
	\begin{align*}
		\gffB_{U,\tau} \overset{\law}{=} \gffB_{U,0},
	\end{align*}
\item 	the driving processes satisfy the system of SDEs~(\ref{eq:circular_Dyson}),
$\xi = \xi_{N}\coloneqq\frac{N+2}{\sqrt{\kappa}}-\frac{\sqrt{\kappa}}{2}$,
$\chi=\frac{2}{\sqrt{\kappa}}-\frac{\sqrt{\kappa}}{2}$,
and $\zeta = \frac{1}{\sqrt{\kappa}}$.
\end{enumerate}
\end{thm}

\begin{rem}
Theorem~\ref{thm:coupling_multi_curve} is, in fact, a generalization of Theorem~\ref{thm:coupling_single_curve} in two ways.
For the one, Theorem~\ref{thm:coupling_multi_curve} covers the case that $N\geq 1$, but it also says that the coupling uniquely fixes driving processes and the parameters $\xi_{N}$ and $\chi$.
\end{rem}

\begin{rem}
\label{rem:collision_Dyson}
The circular $\beta$-Dyson Brownian motions almost surely do not collide when $\beta\geq 1$, and almost surely collide when $0<\beta<1$.
So, the collision time $T$ is, in fact, almost surely infinite for the former case,
and finite for the latter case.
\end{rem}

\begin{proof}
We only prove that (1) implies (2). Once we understand it, the other direction goes by straightforward computation and similar arguments as the single-curve case.

Let us first pick an open set $U\subset \bD$ and a stopping time $\tau$ such that $\bP [\tau\leq \tau_{U}]=1$.
By the domain Markov property of GFF, conditioned over $K_{\tau} = \bD \setminus \bD_{\tau}$, we know that the harmonic function $\harmMAt{N}{\tau}$ is extracted as 
\begin{align*}
	\harmMAt{N}{\tau} = \bE [\gffB_{U,0}|\cF_{K_{\tau}^{+}}].
\end{align*}
Let us take test functions $f_{1},f_{2}\in C^{\infty}_{0}(\bD)$ such that $\supp(f_{i})\subset U$, $i=1,2$.
Then, we know that $(h_{U,0},f_{i})$, $i=1,2$ are Gaussian variables of covariance
\begin{align*}
	\int_{\bD\times \bD}f_{1}(z)G_{\bD}(z,w)f_{2}(w)\leb(dz)\leb(dw).
\end{align*}
By the coupling, sampling them is equivalent to sampling $(\gffD_{\bD_{\tau}},f_{i})$, $i=1,2$ and add to them $(\harmMAt{N}{\tau},f_{i})$, $i=1,2$ that are, due to the domain Markov property, conditionally independent of $\gffD_{\bD_{\tau}}$ given $K_{\tau}$.
Therefore, the covariance of the Gaussian variables $(\harmMAt{N}{\tau},f_{i})$, $i=1,2$ is given by 
\begin{align*}
	\int_{\bD\times \bD}f_{1}(z)(G_{\bD}(z,w) - G_{\bD_{\tau}}(z,w))f_{2}(w)\leb(dz)\leb(dw).
\end{align*}

Now, for $z\in \bD$, we take a sequence of open sets $U_{i}$, $i\in \bN$ such that $\bigcap_{i\in\bN}U_{i} = \{z\}$.
For each $i\in \bN$, the process $(\harmMAt{N}{t\wedge \tau_{U_{i}}}(z):t\ge 0)$ is a martingale.
Since $\tau_{U_{i}} \to \tau_{z}$ as $i\to\infty$, we have that $(\harmMAt{N}{t\wedge \tau_{z}}(z):t\geq 0)$ is a local martingale.
Similarly, for $z,w\in \bD$, we take a sequence of open sets $U_{i}$, $i\in\bN$ such that $\bigcap_{i\in\bN}U_{i}=\{z,w\}$ to see that 
\begin{align*}
	\braket{\harmMAt{N}{t\wedge \tau_{U_{i}}}(z),\harmMAt{N}{t\wedge\tau_{U_{i}}}(w)} = G_{\bD}(z,w) - G_{\bD_{t\wedge \tau_{U_{i}}}}(z,w).
\end{align*} 
Note that the right-hand side makes sense at $z=w$ though each Green's function does not.
At the limit $i\to\infty$, we get
\begin{align*}
	d\braket{\harmMAt{N}{t}(z),\harmMAt{N}{t}(w)} = -dG_{\bD_{t}}(z,w),\quad 0\leq t\leq \tau_{z}\wedge \tau_{w}.
\end{align*}

We have just observed that $(\harmMAt{N}{t}(z):t\in [0,\tau_{z}])$ is a local martingale for any $z\in \bD$.
We take generic points $z_{1},\dots, z_{N}$ and a sequence of stopping times $0=\tau_{0}<\tau_{1}<\cdots $ such that, in each time interval $[\tau_{n},\tau_{n+1}]$, we can almost surely write each $\Theta^{(i)}_{t}$ as a smooth function of $\harmMAt{N}{t}(z_{j})$, $g_{t}(z_{j})$, and $g'_{t}(z_{j})$, $j=1,\dots, N$ due to the implicit function theorem.
This implies that each $(\Theta^{(i)}_{t}:t\geq 0)$ is a semi-martingale admitting a decomposition
\begin{equation}
\label{eq:semimartingale_decomposition}
	\Theta^{(i)}_{t} = \sfM^{(i)}_{t} + \sfU^{(i)}_{t},
\end{equation}
where $\sfM^{(i)}_{t}$ is the martingale and $\sfU^{(i)}_{t}$ is the drift part of $\Theta^{(i)}_{t}$.

We set
\begin{align*}
	\cpxHarmMAt{N}{t}(z)=& -\frac{2}{\sqrt{\kappa}}\sum_{i=1}^{N}\log (g_{t}(z)-X^{(i)}_{t})+\xi\log (g_{t}(z))-\chi \log (g'_{t}(z)) + \bbmi\zeta\sum_{i=1}^{N}\Theta^{(i)}_{t}
\end{align*}
so that $\harmMAt{N}{t}(z)= \rmIm\, \cpxHarmMAt{N}{t}(z)$.
Since $\cpxHarmMAt{N}{t}(z)$ is holomorphic in $z$, its real and imaginary parts are related by the Cauchy--Riemann identities.
Therefore, the drift term of $\cpxHarmMAt{N}{t}(z)$ must be purely real.

Using the decomposition (\ref{eq:semimartingale_decomposition}) we can compute the increment of $\cpxHarmMAt{N}{t}(z)$ by Ito's formula. Computation is very similar to that in the proof of Theorem~\ref{thm:coupling_single_curve} except that we do not know explicit forms of $\sfM^{(i)}_{t}$ and $\sfU^{(i)}_{t}$ at this point.
Let us start with each of the driving processes $(X^{(i)}_{t}:t\geq 0)$, $i=1,\dots, N$:
\begin{align*}
	dX^{(i)}_{t} = \bbmi X^{(i)}_{t} d\sfM^{(i)}_{t} + \biggl(\bbmi \frac{d\sfU^{(i)}_{t}}{dt} - \frac{1}{2}\frac{d\braket{\sfM^{(i)}_{t},\sfM^{(i)}_{t}}}{dt}\biggr)dt.
\end{align*}

The Loewner equation reads
\begin{equation*}
    \frac{d}{dt}g_{t}(z)=\sum_{i=1}^{N}\Psi (g_{t}(z),X^{(i)}_{t})
\end{equation*}
with the same function $\Psi (z,x)$ in (\ref{eq:radial_Loewner_vect_field}).
Thus, we can first compute, for each $i=1,\dots, N$,
\begin{align}
\label{eq:d_of_diff_i}
	d\log (g_{t}(z) - X^{(i)}_{t}) 
	=&\, \frac{d(g_{t}(z) - X^{(i)}_{t})}{g_{t}(z) - X^{(i)}_{t}} + \frac{(X^{(i)}_{t})^{2}}{2(g_{t}(z) - X^{(i)}_{t})^{2}}d\braket{\sfM^{(i)}_{t},\sfM^{(i)}_{t}} \\
	=&\, \frac{-\bbmi X^{(i)}_{t}}{g_{t}(z)-X^{(i)}_{t}}d\sfM^{(i)}_{t} + K^{(i)}_{t}(g_{t}(z),X^{(i)}_{t})dt, \notag
\end{align}
where 
\begin{align*}
K^{(i)}_{t}(z,\bm{x})=&\, \frac{1}{z-x_{i}}\sum_{j=1}^{N}\Psi (z,x_{j})-\frac{x_{i}}{z-x_{i}}\left(\bbmi \frac{d\sfU^{(i)}}{dt}-\frac{1}{2}\frac{d\braket{\sfM^{(i)}_{t},\sfM^{(i)}_{t}}}{dt}\right)\\
+&\, \frac{x_{i}^{2}}{2(z-x_{i})^{2}}\frac{d\braket{\sfM^{(i)}_{t},\sfM^{(i)}_{t}}}{dt}.
\end{align*}
We prepare here a small technical lemma.
\begin{lem}
We have the following identity of rational functions of $z$:
\begin{equation*}
    \sum_{i,j=1}^{N}\frac{\Psi (z,x_{j})}{z-x_{i}}=\sum_{i=1}^{N}\frac{-2x_{i}^{2}}{(z-x_{i})^{2}}+\sum_{i=1}^{N}\frac{-x_{i}}{z-x_{i}}\Biggl(2\sum_{\substack{j=1 \\j\neq i}}^{N}\frac{x_{i}+x_{j}}{x_{i}-x_{j}}+N+2\Biggr)-N^{2}.
\end{equation*}
\end{lem}
\begin{proof}
The diagonal sum immediately becomes
\begin{equation}
\label{eq:tech_lem_diagonal_sum}
    \sum_{i=1}^{N}\frac{\Psi (z,x_{i})}{z-x_{i}}=\sum_{i=1}^{N}\left(-\frac{2x_{i}^{2}}{(z-x_{i})^{2}}-\frac{3x_{i}}{z-x_{i}}\right)-N
\end{equation}
due to the Laurent expansion (\ref{eq:radial_Loewner_vect_field_Laurent}).
The off-diagonal sum
\begin{equation*}
    -\sum_{\substack{i,j=1 \\ i\neq j}}^{N}\frac{z}{z-x_{i}}\frac{z+x_{j}}{z-x_{j}}
\end{equation*}
converges to $-N(N-1)$ as $z\to\infty$ and has simple poles at $z=x_{i}$, $i=1,\dots, N$.
Hence, it admits the partial faction decomposition:
\begin{align}
\label{eq:tech_lem_off-diagonal_sum}
    -\sum_{\substack{i,j=1 \\ i\neq j}}^{N}\frac{z}{z-x_{i}}\frac{z+x_{j}}{z-x_{j}}
    =&\, \sum_{i=1}^{N}\frac{-x_{i}}{z-x_{i}}\Biggl(\sum_{\substack{j=1 \\j\neq i}}\frac{x_{i}+x_{j}}{x_{i}-x_{j}}+\sum_{\substack{j=1 \\j\neq i}}\frac{2x_{i}}{x_{i}-x_{j}}\Biggr)-N(N-1) \\
    =&\,\sum_{i=1}^{N}\frac{-x_{i}}{z-x_{i}}\Biggl(2\sum_{\substack{j=1 \\j\neq i}}\frac{x_{i}+x_{j}}{x_{i}-x_{j}}+N-1\Biggr)-N(N-1). \notag
\end{align}
Combining (\ref{eq:tech_lem_diagonal_sum}) and (\ref{eq:tech_lem_off-diagonal_sum}), we obtain the desired identity.
\end{proof}

We then sum up (\ref{eq:d_of_diff_i}) over $i=1,\dots, N$ to have 
\begin{align*}
	\sum_{i=1}^{N} d\log (g_{t}(z)-X^{(i)}_{t})
	=&\sum_{i=1}^{N}\frac{-\bbmi X^{(i)}_{t}}{g_{t}(z)-X^{(i)}_{t}} d\sfM^{(i)}_{t} + K_{N,t}(g_{t}(z),\bm{X}_{t}),
\end{align*}
where
\begin{align*}
 K_{N,t}(z,\bm{x}) =&\, \sum_{i=1}^{N}K^{(i)}_{t}(z,x_{i}) \\
	=&\, \sum_{i=1}^{N}\Bigg\{ \frac{x_{i}^{2}}{(z-x_{i})^{2}}\left(-2+\frac{1}{2}\frac{d\braket{\sfM^{(i)}_{t},\sfM^{(i)}_{t}}}{dt}\right) \\
	+&\, \frac{-x_{i}}{z-x_{i}}\Bigg(\bbmi \frac{d\sfU^{(i)}_{t}}{dt}-\frac{1}{2}\frac{d\braket{\sfM^{(i)}_{t},\sfM^{(i)}_{t}}}{dt}+2\sum_{\substack{j=1\\ j\neq i}}^{N}\frac{x_{i}+x_{j}}{x_{i}-x_{j}}+N+2\Bigg)\Bigg\}\\
	-&\, N^{2}.
\end{align*}

Both $\log (g_{t}(z))$ and $\log (g'_{t}(z))$ are easier to handle as they only have drift terms. In fact, we obtain
\begin{equation*}
    \frac{d}{dt}\log (g_{t}(z))=L_{N}(g_{t}(z),\bm{X}_{t})
\end{equation*}
with
\begin{equation*}
    L_{N}(z,\bm{x})=\sum_{i=1}^{N}\frac{1}{z}\Psi (z,x_{i})=\sum_{i=1}^{N}\frac{-2x_{i}}{z-x_{i}}-N
\end{equation*}
and
\begin{equation*}
    \frac{d}{dt}\log (g'_{t}(z))=M_{N}(g_{t}(z),\bm{X}_{t})
\end{equation*}
with
\begin{equation*}
    M_{N}(z,\bm{x})=\sum_{i=1}^{N}\frac{\dee \Psi}{\dee z}(z,x_{i})=\sum_{i=1}^{N}\frac{2x_{i}^{2}}{(z-x_{i})^{2}}-N.
\end{equation*}

Combining the above computation altogether, we get
\begin{align*}
	&\, d\cpxHarmMAt{N}{t}(z) \\
	= &\, \frac{\bbmi}{\sqrt{\kappa}}\sum_{i=1}^{N}\left(\frac{2X^{(i)}_{t}}{g_{t}(z)-X^{(i)}_{t}} + \sqrt{\kappa}\zeta\right)d\sfM^{(i)}_{t} \\
	+&\,\left(-\frac{2}{\sqrt{\kappa}}K_{N,t}(g_{t}(z),\bm{X}_{t})+\xi L_{N}(g_{t}(z),\bm{X}_{t})-\chi M_{N}(g_{t}(z),\bm{X}_{t})+\frac{\bbmi}{\sqrt{\kappa}}\sum_{i=1}^{N}\frac{d\sfU^{(i)}_{t}}{dt}\right)dt.
\end{align*}

Let us look at the pole of order 2 at $z=x_{i}$ in the drift term of $d\cpxHarmMAt{N}{t}(z)$.
Since the residue must vanish, we have 
\begin{align*}
	-\frac{2}{\sqrt{\kappa}}\left(-2+\frac{1}{2}\frac{d\braket{\sfM^{(i)}_{t},\sfM^{(i)}_{t}}}{dt}\right)-2\chi = 0.
\end{align*}
By introducing $\lambda=\frac{2}{\sqrt{\kappa}}(\frac{2}{\sqrt{\kappa}}-\frac{\sqrt{\kappa}}{2}-\chi)$, we can solve it as
\begin{equation*}
	\braket{\sfM^{(i)}_{t},\sfM^{(i)}_{t}}=\kappa (1+\lambda) t,\quad t\geq 0.
\end{equation*}
Using this in the first order pole at $z=x_{i}$, we can see that we must have
\begin{align*}
	\frac{2}{\sqrt{\kappa}}\left(\bbmi \frac{d\sfU^{(i)}_{t}}{dt}-\frac{\kappa}{2}(1+\lambda)+2\sum_{\substack{j=1\\ j\neq i}}^{N}\frac{x_{i}+x_{j}}{x_{i}-x_{j}}+N+2\right) - 2 \xi = 0.
\end{align*}
Since $x_{i}\in\bS$, $i=1,\dots, N$, the fraction $\frac{x_{i}+x_{j}}{x_{i}-x_{j}}$ is purely imaginary for all pairs $i\neq j$.
Thus, the imaginary part of the above equation gives 
\begin{equation*}
	\frac{d\sfU^{(i)}_{t}}{dt}=F^{(i)}(\bm{X}_{t}),\quad F^{(i)}(\bm{x})\coloneqq 2\bbmi \sum_{\substack{j=1 \\j\neq i}}^{N}\frac{x_{i}+x_{j}}{x_{i}-x_{j}},\quad t\geq 0
\end{equation*}
and the real part 
\begin{align*}
	\xi = \frac{N+2}{\sqrt{\kappa}} - \frac{\sqrt{\kappa}}{2}(1+\lambda).
\end{align*}

Under this choice of the drift terms $\sfU^{(i)}_{t}$, $i=1,\dots, N$, we have $\sum_{i=1}^{N}\frac{d\sfU^{(i)}_{t}}{dt} = 0$.
Hence, the increment of $\harmMAt{N}{t}(z)$ reads
\begin{align*}
	d\harmMAt{N}{t}(z)=\frac{1}{\sqrt{\kappa}}\sum_{i=1}^{N}\rmRe \left(\frac{2X^{(i)}_{t}}{g_{t}(z)-X^{(i)}_{t}}+\sqrt{\kappa}\zeta\right)d\sfM^{(i)}_{t}.
\end{align*}

The cross variation between $\harmMAt{N}{t}(z)$ and $\harmMAt{N}{t}(w)$ is thus computed as
\begin{align*}
	&\,d\braket{\harmMAt{N}{t}(z),\harmMAt{N}{t}(w)}\\
	=&\,(1+\lambda) \sum_{i=1}^{N}\rmRe \left(\frac{g_{t}(z)+X^{(i)}_{t}}{g_{t}(z)-X^{(i)}_{t}}+(\sqrt{\kappa}\zeta - 1)\right)\rmRe \left(\frac{g_{t}(w)+X^{(i)}_{t}}{g_{t}(w)-X^{(i)}_{t}}+ (\sqrt{\kappa}\zeta - 1)\right)dt \\
	+&\,\frac{1}{\kappa}\sum_{i\neq j}\rmRe \left(\frac{2X^{(i)}_{t}}{g_{t}(z)-X^{(i)}_{t}}+\sqrt{\kappa}\zeta\right)\rmRe \left(\frac{2X^{(j)}_{t}}{g_{t}(w)-X^{(j)}_{t}}+\sqrt{\kappa}\zeta\right)d\braket{\sfM^{(i)}_{t},\sfM^{(j)}_{t}}.
\end{align*}
This must coincide with $-dG_{\bD}(g_{t}(z),g_{t}(w))$, which is similarly computed as in the single-curve case as
\begin{align*}
	-dG_{\bD}(g_{t}(z),g_{t}(w)) = \sum_{i=1}^{N}\rmRe \left(\frac{g_{t}(z)+X^{(i)}_{t}}{g_{t}(z)-X^{(i)}_{t}}\right)\rmRe \left(\frac{g_{t}(w)+X^{(i)}_{t}}{g_{t}(w)-X^{(i)}_{t}}\right)dt
\end{align*}
Therefore, we get $\lambda = 0$ translating to $\chi = \frac{2}{\sqrt{\kappa}}-\frac{\sqrt{\kappa}}{2}$, $\zeta = \frac{1}{\sqrt{\kappa}}$
and $d\braket{\sfM^{(i)}_{t},\sfM^{(j)}_{t}}=0$ for all pairs $i\neq j$.

To conclude, there exist independent standard Brownian motions $(B^{(i)}_{t}:t\geq 0)$, $i=1,\dots, N$
and the processes $(\Theta^{(i)}_{t}:t\geq 0)$, $i=1,\dots, N$ satisfy the system of SDEs (\ref{eq:circular_Dyson}).
The parameters $\xi$, $\chi$, and $\zeta$ are determined in terms of $\kappa$ as desired.
\end{proof}

\subsection{Classical limit $\kappa \to 0$}
\label{sect:coupling_classical}
In Theorem~\ref{thm:coupling_multi_curve}, we assumed $\kappa>0$.
Let us briefly discuss the limit $\kappa \to 0$, which we would like to call the classical limit because it corresponds to the classical limit of CFT at large central charge.
Assume that the equivalent conditions in Theorem~\ref{thm:coupling_multi_curve} are satisfied.
In the limit $\kappa \to 0$, the circular $\sqrt{8/\kappa}$-Dyson Brownian motions converge to a deterministic particle system $(X^{(i)}_{t}=e^{\bbmi \Theta^{(i)}_{t}}:t\geq 0)$, $i=1,\dots, N$ satisfying 
\begin{align}
\label{eq:calogero-sutherland}
	\frac{d}{dt}\Theta^{(i)}_{t} = 2\sum_{\substack{j=1\\j\neq i}}^{N}\cot\left(\frac{\Theta^{(i)}_{t}-\Theta^{(j)}_{t}}{2}\right),\quad t\geq 0,\quad i=1,\dots, N.
\end{align}
In consequence, the Loewner chain $(\SLEgAt{t}:t\geq 0)$ and the stopping time $\tau_{U}$ for any open $U\subset \bD$ are deterministic.

Note that the particle system obeying \eqref{eq:calogero-sutherland} is the Calogero--Sutherland system with the Hamiltonian given by
\begin{align*}
	\cH (\bm{p},\bm{\theta}) = \sum_{i=1}^{N}\frac{p_{i}^{2}}{2} - \sum_{1\leq i < j \leq N} \frac{4}{\sin^{2}\left(\frac{\theta_{i}-\theta_{j}}{2}\right)}.
\end{align*}
In fact, when we set $P_{t}^{(i)}\coloneqq \frac{d\Theta^{(i)}_{t}}{dt}$, $t\geq 0$, $i=1,\dots, N$, we have the closed system of equations
\begin{align*}
	\frac{d}{dt}P^{(i)}_{t} = -\frac{\dee \cH}{\dee \theta_{i}}(\bm{P}_{t},\bm{\Theta}_{t}),\quad \frac{d}{dt}\Theta^{(i)}_{t} = \frac{\dee \cH}{\dee p_{i}}(\bm{P}_{t},\bm{\Theta}_{t}),\quad t\geq 0,\quad i=1,\dots, N.
\end{align*}

As for the GFF, for each $U\subset \bD$ and $t\in [0,\tau_{U}]$, we get
\begin{align*}
\lim_{\kappa \to 0}\sqrt{\kappa}\,\gffB_{U,t} &\,= \lim_{\kappa \to 0}\sqrt{\kappa}\, \harmMAt{N}{t}\\
 &\,= -2\sum_{i=1}^{N}\arg (\SLEgAt{t}(\cdot)-X^{(i)}_{t}) + (N+2)\arg \SLEgAt{t}(\cdot) - 2 \arg \SLEgAt{t}'(\cdot) \\
 &\,\eqqcolon\harmClMAt{N}{t}.
\end{align*}
Due to (1) of Theorem~\ref{thm:coupling_multi_curve}, $\harmClMAt{N}{t}(z)$ is constant in $t\in [0,\tau_{U}]$ for all $z\in U$.
This can be also directly seen from the fact that $(\harmMAt{N}{t}:t\in [0,\tau_{U}])$ is a local martingale (see the proof of Theorem~\ref{thm:coupling_multi_curve}).
In other words, $(\harmClMAt{N}{t}(z):t\in [0,\tau_{U}])_{z\in U}$ is a family of {\it integrals of motion} for the Loewner chain $(\SLEgAt{t}:t\geq 0)$.

This is an example of integrals of motion that arise as the classical limit of martingale observables originated from CFT.
The integrals of motion of this kind have been studied extensively in~\cite{alberts2024pole,makarov2024multiple}.

\section{Hydrodynamic limit}
\label{sect:hydro}
In the previous section, we saw that coupling between radial multiple SLE and GFF is possible if and only if the radial multiple SLE is driven by the circular Dyson Brownian motions,
and concluded that the circular Dyson Brownian motions are the canonical driving processes for radial multiple SLE.
This result, along with the fact that the circular Dyson Brownian motions are a stochastic log-gas, motivates us to study the hydrodynamic limit of the radial multiple SLE, the law of large numbers as $N$ tends to infinity.

\subsection{Hydrodynamic limit of circular Dyson Brownian motions}
\label{sec:PDE}
We start with the hydrodynamic limit of circular Dyson Brownian motions.
Let $\beta>0$ be fixed. 
For each $N\geq 1$, the $N$-particle circular $\beta$-Dyson Brownian motions $(X^{(i)}_{t}:t\geq 0)$, $i=1,\dots, N$ give rise to a measure-valued process 
\begin{align*}
	\DysonEmpAt{N}{t} = \frac{1}{N}\sum_{i=1}^{N}\diracMeasAt{X^{(i)}_{t/N}},\quad t\geq 0.
\end{align*}
Notice the time-change $t \to t/N$ on the particles.
The following is a consequence of \cite[Theorem 4.1]{cepa2001brownian}.
%%%%%%%%%%%%%%%%%%%%
\begin{prop}
\label{thm:Burgers1}
Assume that $\DysonEmpAt{N}{0}$ converges weakly to $\mu_0(\cdot)$ as $N \to \infty$.
The sequence of measure-valued processes $(\DysonEmpAt{N}{t}:t\geq 0)$, $N\geq 1$ almost surely weakly converges to a unique measure-valued process $(\mu_{t}(\cdot):t\geq 0)$ as $N\to\infty$.
Furthermore, the limit process is characterized by its 
circular Stieltjes transform (Helgrotz transform)
\begin{equation}
M_t(z) := \int_{\bS} 
\frac{x+z}{x-z} \mu_t(dx),
\quad
z \in \C \setminus \supp \, \mu_t
\quad t\geq 0, 
\label{eq:Stieltjes}
\end{equation}
that solves the partial differential equation 
\begin{equation}
\frac{\partial M_t(z)}{\partial t}
+2 z M_t(z) \frac{\partial M_t(z)}{\partial z}=0
\label{eq:Burgers}
\end{equation}
under the initial condition
\begin{equation}
M_0(z) = \int_{\bS} 
\frac{x+z}{x-z} \mu_0(dx), 
\quad
z \in \C \setminus \supp \, \mu_0. 
\label{eq:Burgers0}
\end{equation}
\end{prop}
\vskip 0.3cm
%%%%%%%%%%%%%%%%
We may think of \eqref{eq:Burgers} as 
the \textit{complex Burgers equation in the inviscid limit} on the unit circle $\bS$ that has its origin in hydrodynamics.
This is the reason for us calling the limit process $(\mu_t(\cdot) : t \geq 0)$  the \textit{hydrodynamic limit} of the circular Dyson Brownian motions.
Notice that the dependence on the parameter $\beta >0$ disappears in this limit. 
%%%%%%%%%%%%%%%%%%%%%
\begin{rem}
To be precise, the limit process in \cite[Theorem 4.1]{cepa2001brownian} is characterized by the Burgers equation with a diffusion term that we do not have in~\eqref{eq:Burgers}.
This difference arises from the fact that, in \cite{cepa2001brownian}, the Brownian motions scale differently from ours as $N\to \infty$.
However, it is rather straightforward to see in the proof of~\cite[Theorem 4.1]{cepa2001brownian} that we may just drop the diffusion term with our scaling of the Brownian motions.
Detailed analysis with the same scaling as ours can also be found in~\cite{hotta2021limits}.
\end{rem}

\begin{rem}
As we recalled in Remark~\ref{rem:collision_Dyson}, the circular $\beta$-Dyson Brownian motions almost surely collide when $0<\beta<1$ in their original time scale.
Nonetheless, after the time change $t\to t/N$, the limit process $(\mu_{t}(\cdot):t\geq 0)$ can evolve till infinity.
\end{rem}

\begin{lem}
\label{thm:Burgers2}
The solution of \eqref{eq:Burgers} under
the initial condition \eqref{eq:Burgers0} is given by
\[
M_t(z)=M_0(e^{-2 t M_t(z)} z),
\quad z \in \C \setminus \supp \, \mu_t,
\quad t \geq 0.
\]
\end{lem}
%%%%%%%%%%%%%%%%%%%%%%%
%%%%%%%%%%%%%%%%%%%%%%%
\begin{proof}
We apply the method of characteristics.
In our case, we can simply use the time $t$ as the parameter for the characteristic curves.
In fact, when we set $z(t)=we^{2tM_{0}(w)}$, $t\geq 0$ for a fixed $w$,
$M_{t}(z(t))$ stays constant in $t$.
In particular, we have 
\begin{align*}
	M_{t}(z(t)) = M_{0}(w) = M_{0}(e^{-2tM_{t}(z(t))}z(t)),\quad t\geq 0.
\end{align*}
For each $z$, we may find $w$ such that $z(t)=z$.
Therefore, the above equation gives the desired result.
\end{proof}
%%%%%%%%%%%%%%%%%%%%%%%
The limit measure $\mu_t$, $t> 0$ has a smooth density
$\mu_t(dx)=\rho_t(x) dx$ with respect to
the Lebesque measure $dx$ on $\bS$ \cite[Theorem 6.1]{cepa2001brownian}.
For the present circular setting, 
the Sokhotski--Plemelj theorem~\cite[Chapter 7]{ablowitz2003complex} gives
\begin{equation}
\label{eq:stieltjes_inversion}
\rho_t(x)
=\rmRe \left\{
\lim_{r \uparrow 1} \frac{1}{2\pi}
M_t(r x) \right\}
= \rmRe \left\{
\lim_{r \downarrow 1} \frac{1}{2\pi}
M_t(r x) \right\},
\quad x \in \bS, \quad t \geq 0.
\end{equation}

We consider the case \eqref{eq:single_source}. 
That is, the hydrodynamic limit of the
circular Dyson Brownian motions starts from a single source at $z=1$,
\begin{equation}
\mu_0(\cdot)=\diracMeasAt{1}.
\label{eq:single}
\end{equation}
Then, the initial condition \eqref{eq:Burgers0} is given by
\begin{equation}
M_0(z)=\int_{\bS} \frac{x+z}{x-z} \diracMeasAtOf{1}{dx}
=\frac{1+z}{1-z}, \quad z \in \C \setminus\{1\}.
\label{eq:Cauchy0}
\end{equation}
By Lemma \ref{thm:Burgers2}, 
the solution $(M_t(z) : t \geq 0)$
of the inviscid complex Burgers equation \eqref{eq:Burgers} on $\bS$ 
satisfies the functional equation, 
\begin{align*}
M_t(z)= \frac{1+e^{-2 t M_t(z)}z}
{1-e^{-2 t M_t(z)} z}.
\end{align*}
%or equivalently,
%\begin{align*}
%z=e^{2 t M_t(z)}
%\frac{M_t(z)-1}{M_t(z)+1}.
%\label{eq:Cauchy1}
%\end{align*}
The limit process $(\mu_{t}:t\geq 0)$ under this initial condition is known as the multiplicative free Brownian motion~\cite{BercoviciVoiculescu1992, Biane1997free} and was studied in detail by~\cite{Biane1997segal, DemniHmidi2012}.

\begin{prop}
\label{thm:support_Dyson0}
There is a critical time $t_{\mathrm{c}}=1$ such that
\begin{enumerate}
\item
if $0 < t < t_{\mathrm{c}}$, then
the support of $\mu_t$ is an arc given by
\[
\supp \, \mu_t=\{e^{\bbmi \phi} | \,
\phi \in [-\phi_t^{\mathrm{c}}, \phi_t^{\mathrm{c}}] \},
\]
where
\begin{equation}
\phi_t^{\mathrm{c}}=\arccos (1-2t) + 2 \sqrt{t(1-t)},
\label{eq:phitc}
\end{equation}
\item
if $t \geq t_{\mathrm{c}}$, then
the support of $\mu_t$ extends over $\bS$; 
$\supp \, \mu_t= \bS$, where
\begin{enumerate}
\item
at $t=t_{\mathrm{c}}$ 
$\rho_{t_{\mathrm{c}}}(x) > 0$ 
for $x \in \bS \setminus \{-1\}$ and
$\rho_{t_{\mathrm{c}}}(-1)=0$,
\item
at $t \in (t_{\mathrm{c}}, \infty)$, 
$\rho_t(x) > 0$ for all $x \in \bS$.
\end{enumerate}
\end{enumerate}
\end{prop}
%%%%%%%%%%%%%%%%%

%%%%%%%%%%%%%%%%%
\subsection{Hydrodynamic limit of multiple SLE and factorization}
%%%%%%%%%%%%%%%%%%

Associated with the hydrodynamic limit of the circular
Dyson Brownian motions, the following limit theorem
was proved by 
Hotta and Schlei{\ss}inger~\cite{hotta2021limits}.
Let us again fix $\beta>0$.
For $N \in \N$, 
we write the radial multiple SLE driven by the $N$-particle circular $\beta$-Dyson Brownian motions as
$(\mSLEgAt{N}{t} : t \geq 0)$. 

%%%%%%%%%%%%%%%%%%%%%%%%%
\begin{prop}[Hotta and Schlei{\ss}inger~\cite{hotta2021limits}]
\label{thm:hydro_SLE}
Under the same assumption as given in Proposition \ref{thm:Burgers1},
as $N \to \infty$, 
$\mSLEgAt{N}{t/N}$
converges locally uniformly in distribution
to the solution $\mSLEgAt{\infty}{t}$ of the Loewner chain
\begin{equation}
\frac{d}{dt} \mSLEgAt{\infty}{t}(z)
=\mSLEgAt{\infty}{t}(z) M_t(\mSLEgAt{\infty}{t}(z)),
\quad t \geq 0, 
\label{eq:measure_Loewner}
\end{equation}
under the initial condition
$\mSLEgAt{\infty}{0}(z)=z \in \bD$.
\end{prop}
%%%%%%%%%%%%%%%%%%%%%%
\noindent
In contrast to \eqref{eq:radial_multiple_loewner},
Equation~\eqref{eq:measure_Loewner} with
\eqref{eq:Stieltjes} is regarded as
a \textit{measure-driven} radial Loewner equation. 

The following is the version of \cite[Remark 3.11]{delMonacoSchleissinger2016} and \cite[Lemma 2.1]{HottaKatori2018} in the present circular setting.
We write $M_0^{\prime}(y) := \partial M_0(z)/\partial z|_{z=y}$. 
%%%%%%%%%%%%%%%%%%
\begin{lem}
\label{thm:ht}
Define a time-dependent map 
$(h_t : t \geq 0)$ by
\[
h_t(z) = \mSLEgAt{\infty}{t}(z) e^{-2t M_t(\mSLEgAt{\infty}{t}(z))},
\quad t > 0,
\]
with the initial condition
$h_0(z) = \mSLEgAt{\infty}{0}(z)=z \in \bD$.
Then the following equality holds,
\[
M_t(\mSLEgAt{\infty}{t}(z))=M_0(h_t(z)), \quad t \geq 0. 
\]
Moreover, 
$h_t(z)$ solves the equation
\begin{equation}
\frac{d}{d t} h_t(z)
=- \frac{h_t(z) M_0(h_t(z))}
{
1+2 t h_t(z)M_0^{\prime}(h_t(z)) }, \quad t >0. 
\label{eq:h2}
\end{equation}
\end{lem}
\begin{proof}
The first equality is an immediate consequence of Lemma~\ref{thm:Burgers2}.
To derive \eqref{eq:h2}, we compute the time derivative of 
\begin{align}
\label{eq:h_to_g_inf}
	\mSLEgAt{\infty}{t}(z) = h_{t}(z)e^{2tM_{t}(\mSLEgAt{\infty}{t}(z))} = h_{t}(z)e^{2tM_{0}(h_{t}(z))}
\end{align}
in two ways.
By directly differentiating both sides of \eqref{eq:h_to_g_inf}, we get 
\begin{align*}
	\frac{d}{dt}\mSLEgAt{\infty}{t}(z) =&\, \frac{d}{dt}h_{t}(z)\cdot e^{2tM_{0}(h_{t}(z))}
	+ h_{t}(z)\cdot \bigg(2M_{0}(h_{t}(z))+2tM'_{0}(h_{t}(z))\frac{d}{dt}h_{t}(z)\bigg)e^{2tM_{0}(h_{t}(z))} \\
	=&\, \big(1+2th_{t}(z)M'_{0}(h_{t}(z))\big)e^{2tM_{0}(h_{t}(z))}\frac{d}{dt} h_{t}(z) + 2h_{t}(z)M_{0}(h_{t}(z))e^{2tM_{0}(h_{t}(z))}.
\end{align*}
This coincides with \eqref{eq:measure_Loewner} that we can now express as 
\begin{align*}
	\frac{d}{dt}\mSLEgAt{\infty}{t}(z) = h_{t}(z)M_{0}(h_{t}(z))e^{2tM_{0}(h_{t}(z))}.
\end{align*}
Comparing the two expressions, we obtain~\eqref{eq:h2}.
\end{proof}
%%%%%%%%%%%%%%%

\begin{figure}[htbp]
    \begin{tabular}{ccc}
       \centering
        \includegraphics[keepaspectratio, scale=0.3]{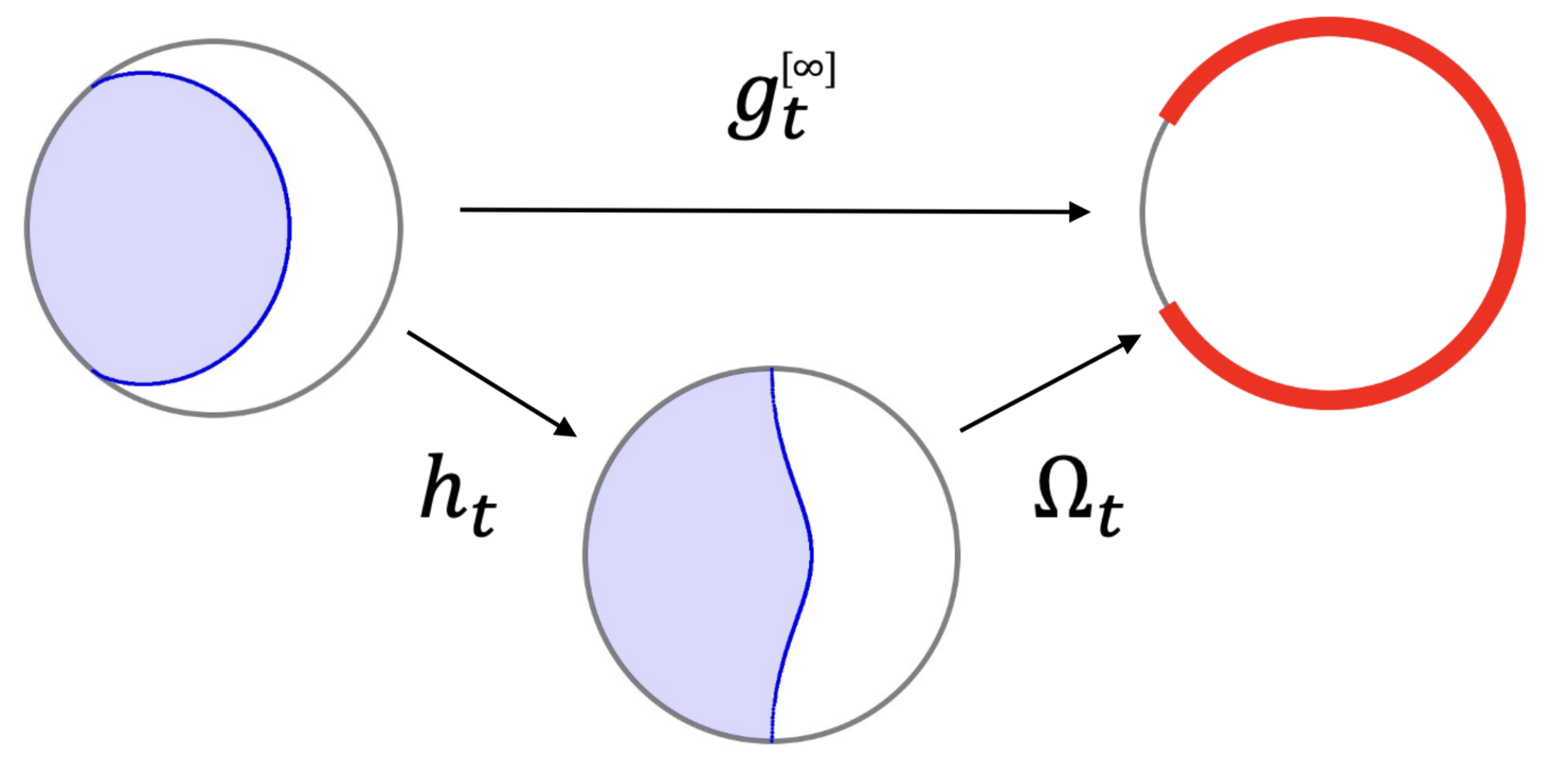}
    \end{tabular}
     \caption{
     In the hydrodynamic limit of the radial multiple SLE, 
     the Loewner chain 
     $\mSLEgAt{\infty}{t}: \mSLEDomAt{\infty}{t} \to \bD$, 
     $t \geq 0$
     is factorized into two steps,
     $h_t: \mSLEDomAt{\infty}{t} \to \widetilde{\bD}_t$ and
     $\Omega_t: \widetilde{\bD}_t \to \bD$.
     Figure shows the system at time 
     $t=0.5 < t_{\mathrm{c}}$ started from $\delta_1$. 
     The support of the hydrodynamic limit $\mu_t$ of the 
     circular Dyson Brownian motion indicated by a red arc
     in the upper right figure is the image of the 
     boundary curve $\gamma_t$ by $\mSLEgAt{\infty}{t}$
     as well as that of $\widetilde{\gamma}_t$ by $\Omega_t$.
     In the text, we use the polar coordinates, 
     $R_t e^{\bbmi \varphi_t} 
     \in \overline{\mSLEDomAt{\infty}{t}}$,
     $r_t e^{\bbmi \theta_t} \in \overline{\widetilde{\bD}_t}$,
     and
     $s_t e^{\bbmi \phi_t} \in \overline{\bD}$, $t \geq 0$,
     respectively. 
     }
     \label{fig:two_step}
  \end{figure}

Equation~\eqref{eq:h_to_g_inf} suggests that the maps $\mSLEgAt{\infty}{t}$, $t\geq 0$ are naturally factorized into (see Figure~\ref{fig:two_step})
\begin{align}
\label{eq:factor_g}
	\mSLEgAt{\infty}{t} = \Omega_{t}\circ h_{t}, \quad t\geq 0
\end{align}
with $\Omega_{t}$ defined by the formula
\begin{equation*}
\Omega_t(z) :=
z e^{2t M_0(z)},
\quad t \geq 0.
\end{equation*}
In what follows, we will study the two maps separately to gain knowledge about the Loewner chains $(\mSLEgAt{\infty}{t}:t\geq 0)$ as well as $(\mSLEDomAt{\infty}{t}:t\geq 0)$.

Let us define $\wtilde{\bD}_{t}\coloneqq h_{t}(\bD_{t})$, $t\geq 0$.
Since $h_{t}$ is univalent on $\bD_{t}$, it is clear that $h_{t}\colon \bD_{t} \to \wtilde{\bD}_{t}$ and $\Omega_{t}\colon \wtilde{\bD}_{t} \to \bD$ are both conformal maps for each $t\geq 0$.
We can also see that the origin $0$ is fixed by these maps.
Therefore, $\wtilde{\bD}_{t}$ is equivalently characterized as the neighbor of the origin that $\Omega_{t}$ conformally maps to $\bD$.

%%%%%%%%%%%%%%%%%%%%%%%%%%%%%%
%%%%%%%%%%%%%%%%%%%%%%%%%%%%%%%

%%%%%%%%%%%%%%%%
\subsection{Solution for the Loewner chain}
\label{sec:Lambert}
%%%%%%%%%%%%%%%%
%%%%%%%%%%%%%%%%%%%%%%%%%%
The goal of this subsection is to provide a proof of Theorem~\ref{thm:hydro_sketch}.
From now on, we consider the case \eqref{eq:single},
that is $(\mu_t: t \geq 0)$  
starts from a single source at $z=1$.
Since \eqref{eq:Cauchy0} gives
$M_0^{\prime}(y)=2/(1-y)^2$, 
\eqref{eq:h2} is written as
\begin{equation}
\frac{d}{d t} h_t(z)
= - \frac{h_t(z) (1-h_t(z))^2}
{1+h_t(z)(h_t(z)+4t-2)}.
\label{eq:h3}
\end{equation}

%%%%%%%%%%%%%%%%%%%%%%%%%%%%%%%%
\begin{lem}
\label{thm:ht3}
Equation \eqref{eq:h3} is solved under
the initial condition $h_0(z)=z \in \bD$ as
\begin{equation}
t=\left( \frac{1-h_t(z)}{1+h_t(z)} \right)^2
\log \left[
\frac{(1-h_t(z))^2 z}{h_t(z) (1-z)^2} \right].
\label{eq:t_h}
\end{equation}
\end{lem}
%%%%%%%%%%%%%%%%%%%
\begin{proof}
Equation \eqref{eq:h3} is written as
\begin{equation}
dt=-\frac{4t}{1-h_t(z)^2} d h_t(z)
- \frac{1-h_t(z)}{h_t(z)(1+h_t(z))} d h_t(z).
\label{eq:eq1}
\end{equation}
We use the constant variation method. 
First we notice that the truncated equation
\[
dt=-\frac{4t}{1-h_t(z)^2} d h_t(z) 
\quad \iff \quad 
\frac{dt}{2t} = -\frac{dh_t(z)}{1-h_t(z)}
- \frac{dh_t(z)}{1+h_t(z)}
\]
is solved as
\begin{equation}
t = c_1 \left( \frac{1-h_t(z)}{1+h_t(z)} \right)^2,
\label{eq:eq2}
\end{equation}
where $c_1$ is a constant. Then we consider the variation 
of $c_1$ and assume $c_1=c_1(h_t(z))$. 
By \eqref{eq:eq1} and \eqref{eq:eq2}, 
we have the relation, 
\[
d c_1=-\frac{1+h_t(z)}{h_t(z)(1-h_t(z)} dh_t(z).
\]
It is readily integrated as
$c_1 = \log \{ c_2 (1-h_t(z))^2/h_t(z) \}$
with a constant $c_2$. The initial condition
$h_0(z)=z$ determines it as
$c_2=z/(1-z)^2$. If we put the obtained $c_1$ into 
\eqref{eq:eq2}, then we obtain \eqref{eq:t_h}.
\end{proof}
%%%%%%%%%%%%%%%%%%%%%%%%%%%%
\noindent
We see that \eqref{eq:t_h} is written as
\begin{align}
& t \left[ 1 + \frac{4 h_t(z)}{(1-h_t(z))^2} \right]
=\log \left[ \frac{(1-h_t(z))^2}{h_t(z)}
\cdot \frac{z}{(1-z)^2} \right]
\nonumber\\
\iff \quad
&\Lambda_t(h_t(z))
e^{\Lambda_t(h_t(z))}
= \Lambda_t(z) e^{-t}, 
\label{eq:t_h2}
\end{align}
where $\Lambda_t(z)$ is defined by \eqref{eq:Lambda}. 

Equation~\eqref{eq:t_h2} can be solved by using the complex Lambert function similarly to~\cite{HottaKatori2018}.
The Lambert $W$-function is defined as the inverse
function of the mapping $x \mapsto x e^x$,
which has two real branches with a branching point 
at $(-e^{-1}, -1)$ in the real plane $(x, W) \in \R^2$
\cite{corless1996lambert,veberivc2012lambert}.
We take the upper branch $W_0(x)$ defined for
$x \in [-e^{-1}, \infty)$; 
\begin{equation}
W_0(x) e^{W_0(x)}=x,
\label{eq:Lambert1}
\end{equation}
where
$W_0(0)=0$, $W_0(e)=1$, and
$W_0(x) \sim x$ as $x \to 0$.
It admits the series expansion
\begin{align*}
W_0(x)=\sum_{n=1}^{\infty} \frac{(-n)^{n-1}}{n!} x^n,
\end{align*}
whose radius of convergence is $e^{-1}$.
This allows us to extend it to the complex function $W_0(z)$, $|z|<e^{-1}$.
In fact, $W_{0}$ is analytically continued to $\bC \setminus (-\infty, -e^{-1}]$ (see~\cite{corless1996lambert} for detail).
Comparing \eqref{eq:t_h2} with \eqref{eq:Lambert1},
we conclude
\[
W_0(\Lambda_t(z) e^{-t}) = \Lambda_t(h_t(z)) = 4t\frac{h_{t}(z)}{(1-h_{t}(z))^{2}}.
\]
By the initial condition $h_0(z)=z$, 
this provides
\begin{align}
\label{eq:h_t_explicit_formula}
h_t(z)=1+
\frac{2t}{W_0(\Lambda_t(z) e^{-t})}
\Big(1-
\sqrt{1+W_0(\Lambda_t(z) e^{-t})/t} \Big),
\quad t \geq 0.
\end{align}

As for the map $\Omega_{t}$, $t\geq 0$, we now have the explicit formula
\begin{equation}
\Omega_t(z)
=z \exp \left( 
2t \frac{1+z}{1-z} \right),
\quad t \geq 0,
\label{eq:Omega1}
\end{equation}
under the initial condition \eqref{eq:single} along with \eqref{eq:Cauchy0}.
Putting \eqref{eq:h_t_explicit_formula} and \eqref{eq:Omega1} together by the factorization~\eqref{eq:factor_g}, we complete the proof of Theorem~\ref{thm:hydro_sketch}.

%%%%%%%%%%%%%%%%
\subsection{Analysis of the map $\Omega_t$}
\label{sec:Omega}
%%%%%%%%%%%%%%%%
It remains to prove Theorem~\ref{thm:hydro_domain}.
As a step towards it, we first obtain an analogous result for $\wtilde{\bD}_{t}$, $t\geq 0$ through the study of the map 
$\Omega_t, t \geq 0$
given by \eqref{eq:Omega1}. 
In the polar coordinate,
$z=r e^{\bbmi \theta}$, 
$r \in [0, 1]$, 
$\theta \in (-\pi, \pi]$,
we have
\begin{equation*}
\Omega_t(r e^{\bbmi \theta})
= s_t(r, \cos \theta) e^{\bbmi \phi_t(r, \theta)},
\end{equation*}
where
\begin{align}
s_t(r, a) &=
r \exp \left[ 
\frac{2t(1-r^2)}{1-2 r a + r^2} \right],
\label{eq:Omega2b}
\\
\phi_t(r, \theta) &=
\theta+
\frac{4t r \sin \theta}{1-2 r \cos \theta + r^2}.
\label{eq:Omega2c}
\end{align}
We study the condition
$\Omega_t(r e^{\bbmi \theta}) \in \bS$, that is, 
solutions of the equation $s_t(r, \cos \theta)=1$, in detail. 
It is obvious that $\Omega_t(0)=0$,
and \eqref{eq:Omega2b} satisfies
$s_t(1,a)=1$ for $a \not=1$.
When $a=1$, we can readily verify that
$s_t(r, 1)$ is increasing in $r$ and 
$s_t(r, 1) \uparrow \infty$ as $r \uparrow 1$.
Otherwise, we have
\[
\frac{\partial s_t(r, a)}{\partial r}
=U_{t,a}(r)q_{t,a}^+(r) q_{t,a}^-(r)
\]
with
\begin{align*}
U_{t,a}(r)&=\frac{1}{(1-2ra+r^2)^2} \exp \left[ \frac{2t(1-r^2)}{1-2ra+a^2} \right], \\
q_{t,a}^{\pm}(r)&=r^2+2 [-a(1-t) \pm \sqrt{A_{t,a}} ] r+1, 
\end{align*}
where
$A_{t, a}=t [2(1-a^2)+a^2 t] \geq 0$ for $t \geq 0$ and 
$a \in [-1, 1)$.
It is obvious that $U_{t,a}(r)>0$, $0\leq r<1$.
By straightforward calculation, we can prove the following.
\begin{enumerate}
\item
For $t > 0$ and $a \in [-1, 1)$, 
$q_t^{+}(r)$ does not have any root in $(0, 1)$.
\item
When $0 < t \leq 1$, if $-1 \leq a < 1-2t$,
$q_t^-(r)$ does not have any root in $(0, 1)$,
but if $1-2t \leq a < 1$,
$q_t^-(r)$ has a unique root
$r_{t, a}^{\ast} \in (0,1)$.
\item
When $t>1$,
$q_t^-(r)$ has a unique root
$r_{t, a}^{\ast} \in (0,1)$
for any $a \in [-1, 1)$.
\end{enumerate}
Let $a \in [-1, 1)$ be fixed.
Since $s_t(0,a)=0 < s_t(1, a)=1$, 
if there is only one root $r_{t,a}^{\ast} \in (0, 1)$,
$s_t(r_{t,a}^{\ast}, a)>1$
is a maximal value for $r \in (0, 1)$, 
and hence $s_t(r, a)$ hits 1 
at some $r \in (0, r_{t,a}^{\ast})$.
Therefore, the equation
$s_t(r, a)=1$ has one or two solutions in $(0, 1]$,
where one of them is the trivial solution $r=1$,
and other solution is in $(0, r_{t, a}^{\ast})$,
if exists.
For $\theta \in (-\pi, \pi]$, we define
$r_t(\theta)$ as the smallest solution of
$s_t(r, \cos \theta)=1$.
Obviously, $r_t(\theta)=r_t(2\pi-\theta)$, 
because $s_t(r, \cos \theta)$ only depends on $\cos \theta$.
We can also see that, when $t \leq 1$, there is a critical angle
$\theta_t^{\mathrm{c}} \in [0, \pi]$ such that
$r_t(\theta) < 1$ 
for $-\theta_t^{\mathrm{c}} < \theta < \theta_t^{\mathrm{c}}$,
and $r_t(\theta)=1$ 
for $\theta_t^{\mathrm{c}} \leq \theta \leq \pi$ or
$-\pi < \theta \leq -\theta_t^{\mathrm{c}}$, 
while when $t >1$, 
$r_t(\theta) <1$ for all $\theta \in (-\pi, \pi]$.

At this point, we may identify $\wtilde{\bD}_{t}$, $t>0$ as 
\begin{align*}
	\widetilde{\bD}_t = \{r e^{\bbmi \theta} | \, 
0 \leq r < r_t(\theta),
\theta \in (-\pi, \pi] \},
\end{align*}
but we would be interested in the part of $\partial \wtilde{\bD}_{t}$ that runs in $\bD$.
The key observation is that, if $r_{t}(\theta)<1$,
the equation $s_t(r_{r}(\theta), \cos \theta)=1$ with \eqref{eq:Omega2b} 
is written as
\[
\theta = \arccos F_t(r_{t}(\theta)),
\]
where
\begin{equation}
F_t(r)=\frac{t(1-r^2)}{r \log r} + \frac{1+r^2}{2r}, \quad 0 < r < 1.
\label{eq:Ftb}
\end{equation}
Let us set $I_{t}\coloneqq F_{t}^{-1}([-1,1])$.
We can see that $F_{t}(r)$ is monotonically decreasing in $r\in I_{t}$ regardless of $t>0$.
When $0< t\leq 1$, we have $I_{t} = [r_{t}^{\min},1)$ with $r_{t}^{\min}\in (0,1)$ defined as the unique solution of
\begin{align}
F_t(r_t^{\min})=1 \quad \iff \quad
t=-\frac{1-r_t^{\min}}{2(1+r_t^{\min})}
\log r_t^{\min},
\label{eq:rho_min}
\end{align}
and
\[
\lim_{r \uparrow 1}
F_t(r)=1-2t.
\]
Therefore, the critical angle is explicitly determined as
\begin{align}
\theta_t^{\mathrm{c}}=\arccos (1-2t) \in (0, \pi],\quad t\leq 1.
\label{eq:theta_c}
\end{align}
When $t>1$, we have $I_{t}=[r_{t}^{\min}, r_{t}^{\max}]$,
where $r_{t}^{\min}$ is given by the same formula \eqref{eq:rho_min}
and $r_{t}^{\max}\in (0,1)$ is the unique solution of 
\begin{align}
F_t(r_t^{\max})=-1 \quad \iff \quad
t=-\frac{1+r_t^{\max}}{2(1-r_t^{\max})}
\log r_t^{\max}.
\label{eq:rho_max}
\end{align}
%%%%%%%%%%%%%

\begin{figure}[htbp]
    \begin{tabular}{ccc}
    \hskip -2.5cm
      \begin{minipage}[t]{0.63\hsize}
       \centering
        \includegraphics[keepaspectratio, scale=0.22]{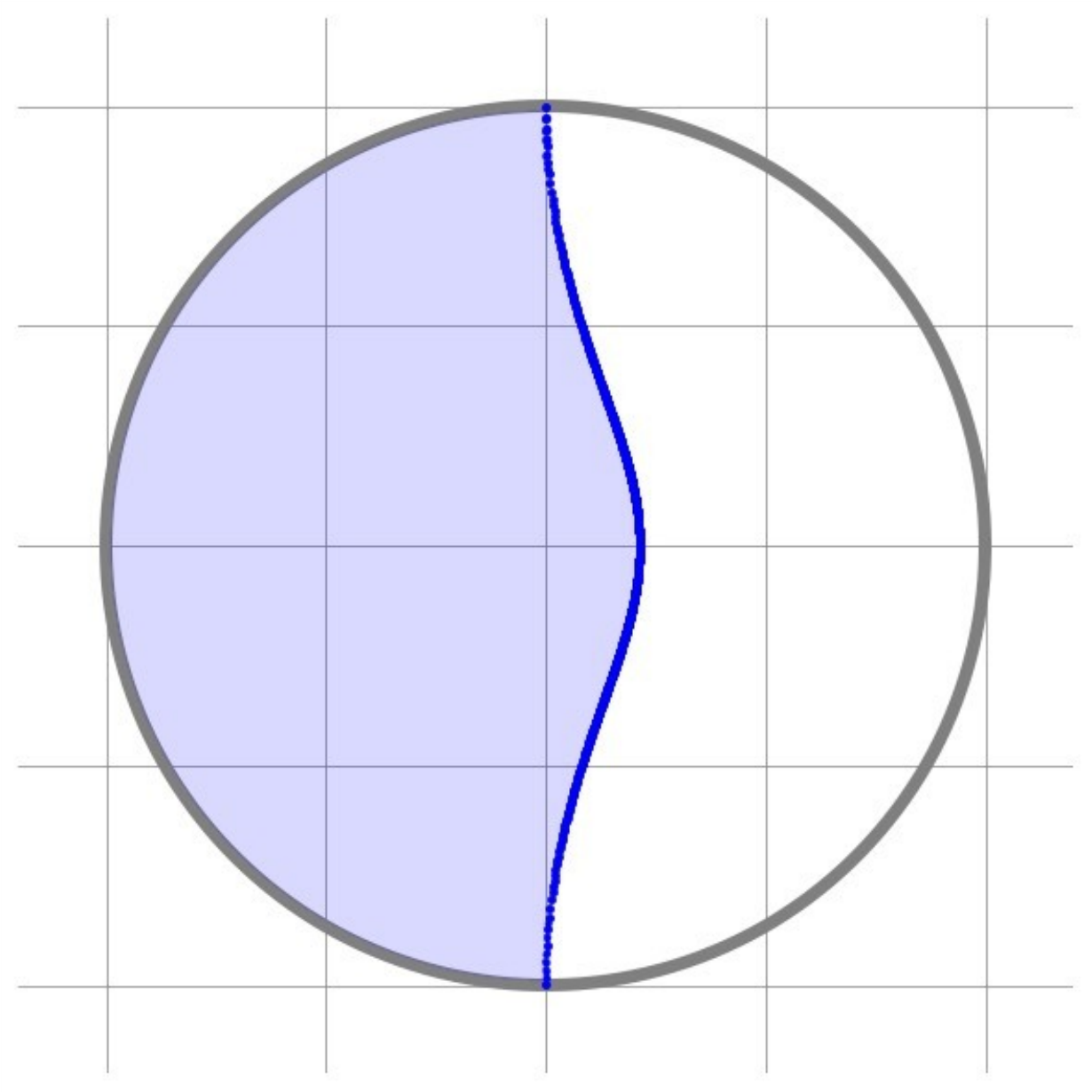}
        \subcaption{}
        \label{fig:SLE_pre_t=0.5}
      \end{minipage} &
    \hskip -5.3cm
      \begin{minipage}[t]{0.63\hsize}
       \centering
        \includegraphics[keepaspectratio, scale=0.22]{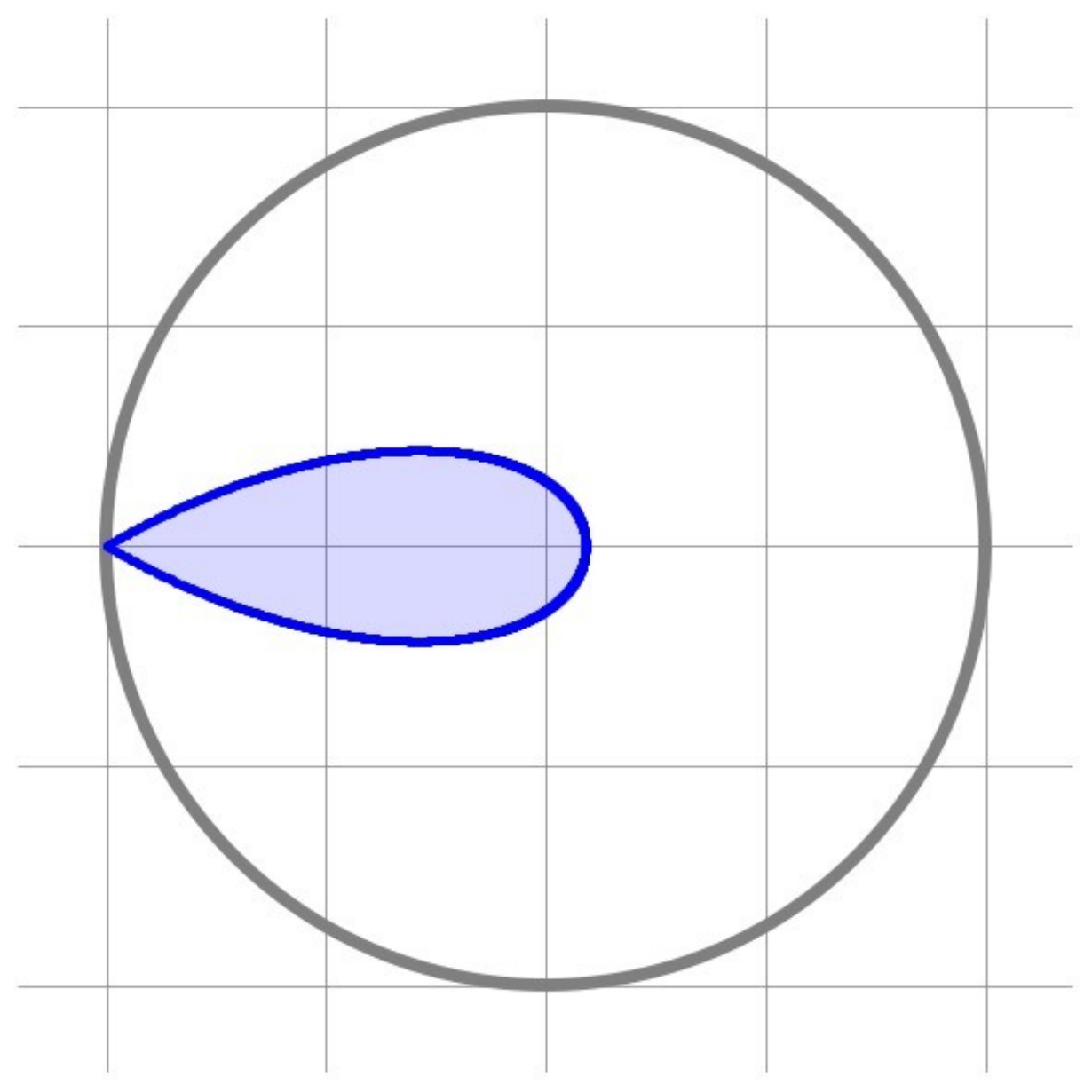}
        \subcaption{}
        \label{fig:SLE_pre_t=1}
      \end{minipage} &
          \hskip -5.3cm
      \begin{minipage}[t]{0.63\hsize}
       \centering
        \includegraphics[keepaspectratio, scale=0.22]{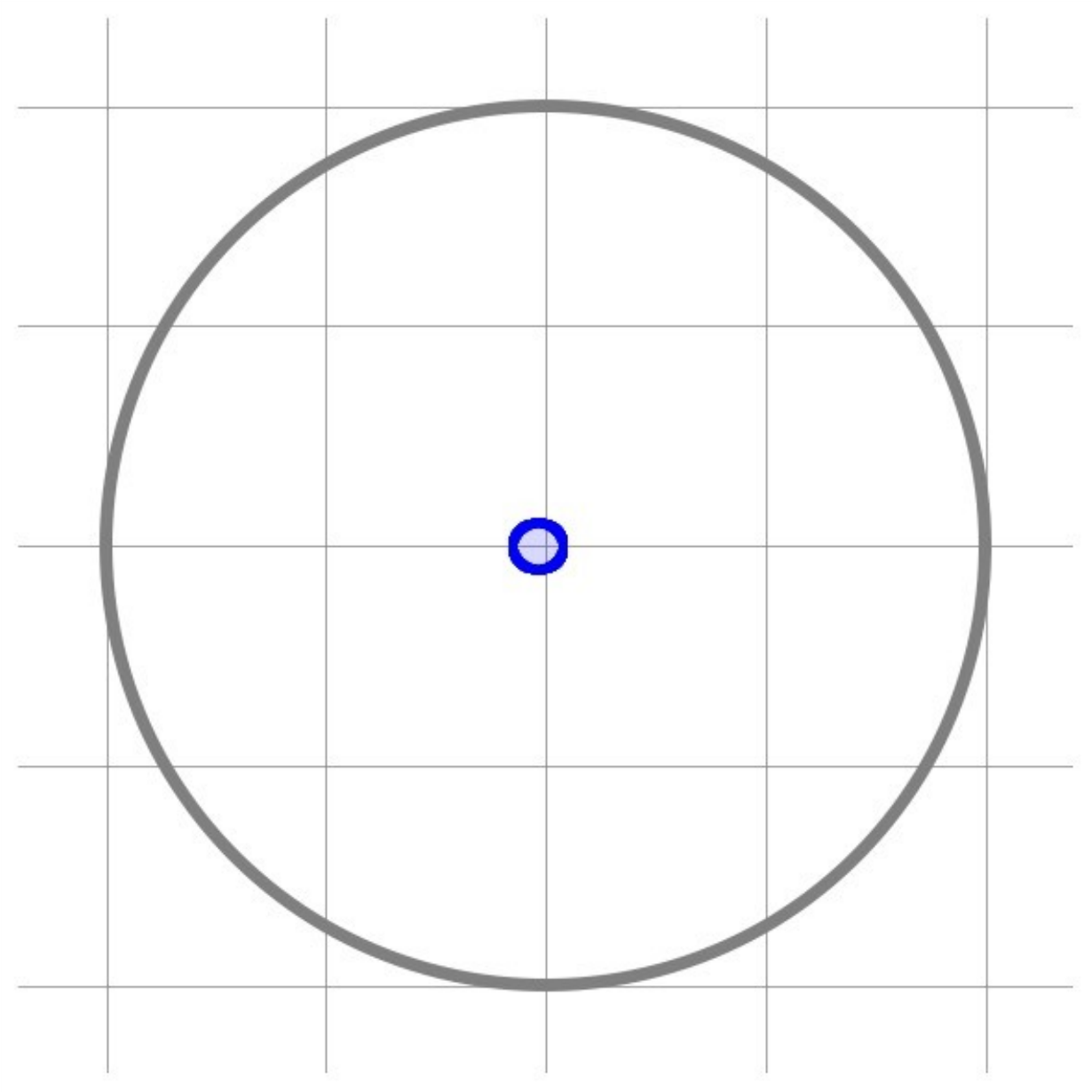}
        \subcaption{}
        \label{fig:SLE_pre_t=1.5}
      \end{minipage} 
    \end{tabular}
     \caption{
     The domains 
     $\widetilde{\bD}_t$ of the map
     $\Omega_t$ are shown 
     by the shaded subdomains of $\bD$
     and the curves $\widetilde{\gamma}_t$ 
     are drawn by tick lines for
     (a) $t=0.5$, 
     (b) $t=t_{\mathrm{c}}=1$ (critical time), 
     (c) $t=1.5$, respectively.
     }
     \label{fig:SLE_pre}
  \end{figure}

With the critical time $t_{\mathrm{c}}=1$, 
we define a simple curve in $\overline{\bD}$ by
\begin{align}
\widetilde{\gamma}_t
&=\begin{cases}
\displaystyle{
\left\{ r e^{\bbmi \theta} \left| \,
\theta = \arccos F_t(r), 
r \in [r_t^{\min}, 1)
\right. \right\} \cup \{e^{\pm \bbmi \theta_{t}^{\mathrm{c}}}\}
}, &
\mbox{if $0 \leq t \leq t_{\mathrm{c}}$},
\cr
\displaystyle{
\left\{ r e^{\bbmi \theta} \left| \,
\theta =\arccos F_t(r), 
r \in [r_t^{\min}, r_t^{\max}]
\right. \right\}
}, &
\mbox{if $t_{\mathrm{c}} < t < \infty$}.
\end{cases}
\label{eq:omega_t2}
\end{align}
Note that, unless $r=r_{t}^{\min}$ or $r_{t}^{\max}$,
there are two $\theta$'s satisfying the equation for a given $r$.
It is clear that $\wtilde{\gamma}_{t}$ splits $\bD$ into two parts.
The analogue of Theorem~\ref{thm:hydro_domain} for $\wtilde{\bD}_{t}$ goes as follows (see Figure \ref{fig:SLE_pre} for illustration).

\begin{prop}
\label{prop:tilde_D_t_characterization}
At each $t\geq 0$, $\wtilde{\bD}_{t}$ is the connected component of $\bD\setminus \wtilde{\gamma}_{t}$ that contains $0$.
\end{prop}

Finally, if we set $r=1$ in \eqref{eq:Omega2c}, we get
\[
\phi_t(1, \theta)
=\theta+\frac{2t \sin \theta}{1-\cos \theta}
=\theta+2 t \cot(\theta/2).
\]
Since \eqref{eq:theta_c} gives
$\sin (\theta_t^{\mathrm{c}}/2)=\sqrt{t}$ and
$\cos (\theta_t^{\mathrm{c}}/2)=\sqrt{1-t}$, 
we have
\begin{equation*}
\phi_t(1, \theta_t^{\mathrm{c}})
= \arccos (1-2t) + 2 \sqrt{t(1-t)}, 
\quad 0 \leq t \leq t_{\mathrm{c}}.
%\label{eq:phi_c}
\end{equation*}
This is exactly equal to 
$\phi_t^{\mathrm{c}}$ given by \eqref{eq:phitc}.
Then we can conclude the following.
%%%%%%%%%%%%%%%%%%%%%
\begin{prop}
\label{thm:support_Dyson}
For $t>0$, 
the support of the hydrodynamic limit 
of the circular Dyson Brownian motions
$\mu_t$ on $\bS$ is the image of
the curve $\widetilde{\gamma}_t$ in $\overline{\bD}$ 
by the map $\Omega_t$; 
\[
\supp \, \mu_t= \Omega_t(\widetilde{\gamma}_t),
\quad t>0.
\]
\end{prop}

%%%%%%%%%%%%%%%%%%%%%%%%%%%%%%
%%%%%%%%%%%%%%%%%%%%%%%%%%%%%%%

%%%%%%%%%%%%%%%%
\subsection{Boundary curve of the SLE hull}
\label{sec:boundary}
%%%%%%%%%%%%%%%%

Let us finalize the proof of Theorem~\ref{thm:hydro_domain}.
We shall define $\what{\gamma}_{t}\coloneqq h_{t}^{-1}(\wtilde{\gamma}_{t})$, $t\geq 0$ and show that it coincides with $\gamma_{t}$, $t\geq 0$ defined in \eqref{eq:SLE_hull}.
Combined with Proposition~\ref{prop:tilde_D_t_characterization}, we may conclude Theorem~\ref{thm:hydro_domain}.

Let us fix $t\geq 0$.
The key identity for determining $\what{\gamma}_{t}$ is \eqref{eq:t_h2}. 
In order to exploit it effectively, let us first observe that 
$\Lambda_{t}(e^{\zeta}) = t/\sinh^{2}(\zeta/2)$.
In this new variable, \eqref{eq:t_h2} gives 
\begin{align}
\label{eq:zeta_from_h}
	\sinh^{2}\frac{\zeta}{2} = \frac{te^{-t}}{\Lambda_{t}(h_{t}(e^{\zeta}))e^{\Lambda_{t}(h_{t}(e^{\zeta}))}}.
\end{align}
Let us investigate $\Lambda_{t}(h_{t}(e^{\zeta}))$ for such $\zeta$ that $h_{t}(e^{\zeta})\in \wtilde{\gamma}_{t}$ in more detail.

In the polar coordinates $z=re^{\bbmi \theta}$, 
\eqref{eq:Lambda} has the following expression:
\[
\Lambda_t(re^{\bbmi \theta})
= \frac{4tr}{1-2 r \cos \theta + r^2}
e^{\bbmi \varpi (r, \cos \theta)}, \quad t \geq 0,
\]
with
\[
\cos \varpi (r, a) = \frac{(1+r^2) a - 2r}
{1-2r a+ r^2}.
\]
When $re^{\bbmi \theta}\in \wtilde{\gamma}_{t}$ with $\theta\in [0,\pi]$,
$\Lambda_{t}(re^{\bbmi\theta})$ is a function only of $r$ because $\theta$ is uniquely determined by~\eqref{eq:omega_t2}.
Let us set 
\begin{equation*}
\widetilde{\Lambda}_t(r):=
\Lambda_t(re^{\bbmi \theta}) 
\Big|_{\theta = \arccos F_{t}(r)\in [0,\pi]}
=\widetilde{u}(r) e^{2 \bbmi \widetilde{\psi}_t(r)}
\end{equation*}
with
\begin{align*}
\widetilde{u}(r) &= \frac{4tr}{1-2rF_{t}(r)+r^{2}} =  -\frac{2r \log r}{1-r^2},\\
\cos ( 2\widetilde{\psi}_t(r))
&= \cos \varpi (r,F_{t}(r)) = - \frac{(1-r^2) \log r}{4 t r} - \frac{1+r^2}{2r}.
\end{align*}
Then, \eqref{eq:zeta_from_h} allows us to express $\zeta\in \bC$ such that $e^{\zeta}\in \what{\gamma}_{t}\cap \ol{\bH}$ as 
\begin{align*}
	\zeta = \wtilde{\Phi}_{t}(r) & \coloneqq 2 \arcsinh \left[-\sqrt{
\frac{t e^{-t}}{\widetilde{\Lambda}_t(r) 
e^{\widetilde{\Lambda}_t(r)}} }\right] \\
&=2\arcsinh \left[- \sqrt{ \frac{t e^{-t}}{\widetilde{u}(r)}}
\exp \left( - \bbmi \widetilde{\psi}_t(r)
-\frac{1}{2} \widetilde{u}(r)
e^{2 \bbmi \widetilde{\psi}(r)} \right)\right]
\end{align*}
with some allowed $r$.
Notice that the function \eqref{eq:Phi} is identified as  $\Phi_{t}(x) = \wtilde{\Phi}_{t}(e^{-x})$.
When we set $x_{t}^{\min}=-\log r_{t}^{\max}$ and $x_{t}^{\max} = -\log r_{t}^{\min}$,
the definitions of $r_{t}^{\min}$ and $r_{t}^{\max}$ in \eqref{eq:rho_min} and \eqref{eq:rho_max} are translated to \eqref{eq:x_min_max}.
We have just seen that 
$\what{\gamma}_{t}\cap \ol{\bH} = \gamma_{t}^{+}$,
but we can conclude that $\what{\gamma}_{t} = \gamma_{t}$
due to the reflection symmetry against the real axis.

%%%%%%%%%%%%%%%%%%%
\subsection{Alternative proof of Proposition~\ref{thm:support_Dyson0}}
%%%%%%%%%%%%%%%%%%%
We proved Proposition~\ref{thm:support_Dyson} by matching the previously known Proposition~\ref{thm:support_Dyson0} and a result on the map $\Omega_{t}$, $t\geq 0$.
However, we can also prove Proposition~\ref{thm:support_Dyson} independently and thereby give an alternative proof of Proposition~\ref{thm:support_Dyson0}.
In fact, the map $\Omega_{t}$, $t\geq 0$ is essentially the so-called $\Sigma$-transform of the multiplicative free Brownian motion~\cite{BercoviciVoiculescu1992, Biane1997free}.

Let $M_{t}(z)$, $t\geq 0$ be the circular Stieltjes transform as in Proposition~\ref{thm:Burgers1} with the single source initial condition.
We form a new function 
\begin{align*}
	\eta_{t}(z)\coloneqq \frac{M_{t}(z)-1}{M_{t}(z)+1},
\end{align*}
which is univalent on the unit disk $\bD$ at each $t\geq 0$.
Using this function, the inversion formula \eqref{eq:stieltjes_inversion} reads 
\begin{align*}
	\rho_{t}(x) = \frac{1}{2\pi}\lim_{r\uparrow 1} \frac{1-|\eta_{t}(rx)|^{2}}{|1-\eta_{t}(rx)|^{2}},\quad x\in \bS,\quad t\geq 0.
\end{align*}
Therefore, we have the following:
\begin{align*}
	\supp\mu_{t}=\Big\{x\in \bS\, \Big|\, \lim_{r\uparrow 1}|\eta_{t}(rx)|<1\Big\}.
\end{align*}

Although the Stieltjes transform $M_{t}(z)$, $t\geq 0$ is characterized by a non-linear partial differential equation, the inverse function of $\eta_{t}$, $t\geq 0$ behaves in a much simpler way. Indeed, we have the following. 

\begin{prop}
At each $t\geq 0$, we have $\Omega_{t} = \eta_{t}^{-1}$.
\end{prop}
\begin{proof}
The partial derivative in $t$ of the functional equation 
\begin{align*}
	\eta_{t}(\eta^{-1}_{t}(z)) = z
\end{align*}
yields
\begin{align}
\label{eq:stationary_M_eta}
	\frac{\dee}{\dee t}M_{t}(\eta_{t}^{-1}(z)) = 0.
\end{align}
At $t = 0$, we have $\eta_{0}(z) = z$, and thus $\eta^{-1}_{0}(z) = z$.
The above differential equation~\eqref{eq:stationary_M_eta} means that, for each $z$, the map $t\mapsto \eta^{-1}_{t}(z)$ is the characteristic curve of the complex Burgers equation \eqref{eq:Burgers} passing through $z$ at $t=0$.
Therefore, we must have (see the proof of Lemma~\ref{thm:Burgers2})
\begin{align*}
	\eta_{t}^{-1}(z) = ze^{2t M_{0}(z)} = \Omega_{t}(z)
\end{align*}
as is desired.
\end{proof}

In this way, Proposition~\ref{thm:support_Dyson} is a simple corollary, and the analysis in Subsection~\ref{sec:Omega} proves Proposition~\ref{thm:support_Dyson0}.
It is worth noting that the map $\Omega_{t}$, $t\geq 0$ naturally appeared in the factorization of the Loewner map as in~\eqref{eq:factor_g}.
In this sense, we can say that analysis of multiple SLE contains analysis of dynamical random matrices at the hydrodynamic limit.

%%%%%%%%%%%%%%%%%
\subsection{Edge asymptotics}
%%%%%%%%%%%%%%%%%%

%%%%%%%%%%%%%%%%%%%%%

As long as $t\leq t_{\mathrm{c}}=1$, the boundary curve $\gamma_{t}$ hits the unit circle at $e^{\pm \bbmi \varphi_{t}^{\mathrm{c}}}$ with the critical angle $\varphi_{t}^{\mathrm{c}}$ given in \eqref{eq:varphi_t_c}.
Let us take a closer look at the behaviour of $\gamma_{t}$ in the vicinity of the edges.
We define $R_{t}(\varphi)\in (0,1]$, $\varphi\in [-\varphi_{t}^{\mathrm{c}}, \varphi_{t}^{\mathrm{c}}]$ by $\gamma_{t} = \{R_{t}(\varphi)e^{\bbmi \varphi}|\varphi\in [-\varphi_{t}^{\mathrm{c}}, \varphi_{t}^{\mathrm{c}}]\}$.

\begin{prop}
\label{thm:vicinity2}
%\begin{description}
\begin{enumerate}
\item
For $0 < t < t_{\mathrm{c}}$,
$R_{t}(\varphi)$ behaves in the vicinity of $\varphi_t^{\mathrm{c}}$ as
\begin{equation}
1- R_t(\varphi) \sim b_t (\varphi_t^{\mathrm{c}}-\varphi)^{\nu}
\quad \mbox{as $\varphi \uparrow \varphi_t^{\mathrm{c}}$ \, with}
\quad \nu= \frac{3}{2},
\label{eq:edgeB1}
\end{equation}
where the coefficient 
\[
b_t=\frac{\sqrt{2}}{3 t^{1/4}(1-t)}
\left( \frac{1-t e^{1-t}}{e^{1-t}} \right)^{1/4}
\]
diverges as
\[
b_t \sim \frac{\sqrt{2}}{3} (1-t)^{-\sigma} \uparrow +\infty \quad
\mbox{as $t \uparrow t_{\mathrm{c}}$ \, with} \quad
\sigma=\frac{3}{4}.
\]
\item
At $t=t_{\mathrm{c}}$, we have
\begin{equation}
1- R_{t_{\mathrm{c}}}(\varphi) \sim \frac{1}{\sqrt{3}} 
(\varphi_{t_{\mathrm{c}}}^{\mathrm{c}}-\varphi)
\quad \mbox{as $\varphi \uparrow 
\varphi_{t_{\mathrm{c}}}^{\mathrm{c}}=\pi$}.
\label{eq:edgeB2}
\end{equation}
\end{enumerate}
\end{prop}
%%%%%%%%%%%%%%%%%
\begin{proof}
For $0 < x \ll 1$, \eqref{eq:ux} gives
\[
u(x)=1 - \frac{x^2}{6} + \frac{7 x^4}{360} + \cO(x^6).
\]
\noindent (1)
If $0 < t < t_{\mathrm{c}}$, \eqref{eq:psi} gives
\[
\psi_t(x)=\frac{\pi}{2}
-\frac{1}{2} \left( \frac{1-t}{t} \right)^{1/2} x
-\frac{1}{48 t^{3/2}(1-t)^{1/2}} x^3
+\cO(x^5),
\]
and then \eqref{eq:Phi} is written as
\[
\Phi_t(x)
= 2\arcsinh\left[\bbmi t^{1/2} e^{(1-t)/2}
- c_t^{\mathrm{R}} x^3 - \bbmi c_t^{\mathrm{I}} x^2
+\cO(x^4)\right],
\]
where
\[
c_t^{\mathrm{R}}=\frac{(1-t)^{1/2}}{12 t} e^{(1-t)/2}, \quad
c_t^{\mathrm{I}}=\frac{t^{1/2} (1-t)}{4t} e^{(1-t)/2}.
\]
For the moment, let us put $R_{t}(x)e^{\bbmi \varphi_{t}(x)}\coloneqq e^{\Phi_t(x)}$ so that both $R_{t}(x)$ and $\varphi_{t}(x)$ are regarded as functions of $x$.
Then, \eqref{eq:gamma_t} with \eqref{eq:Phi} 
gives
\begin{align*}
&\, \bbmi t^{1/2} e^{(1-t)/2}
- c_t^{\mathrm{R}} x^3 - \bbmi c_t^{\mathrm{I}} x^2
+\cO(x^4) \\
=&\,\bbmi \sin(\varphi_t^{\mathrm{c}}/2)
-\frac{1}{2} (1-R_t(x)) \cos(\varphi_t^{\mathrm{c}}/2)
- \bbmi \frac{1}{2}(\varphi_t^{\mathrm{c}}-\varphi_t(x))
\cos(\varphi_t^{\mathrm{c}}/2).
\end{align*}
Since we have $\sin(\varphi_t^{\mathrm{c}}/2) = t^{1/2} e^{(1-t)/2}$ from \eqref{eq:varphi_t_c},
we get
\begin{align*}
&1-R_t(x) \sim
\frac{2 c_t^{\mathrm{R}}}{\cos (\varphi_t^{\mathrm{c}}/2)} x^3,
\quad
\varphi_t^{\mathrm{c}}-\varphi_t(x)
\sim \frac{2 c_t^{\mathrm{I}}}{\cos (\varphi_t^{\mathrm{c}}/2)} x^2
\end{align*}
by comparing the real and imaginary parts.
The assertion (1) follows because $\varphi_{t} (x) \uparrow \varphi_{t}^{\mathrm{c}}$ as $x\downarrow 0$.

\noindent (2) 
If $t=t_{\mathrm{c}}$, \eqref{eq:psi} gives
\[
\psi_{t_{\mathrm{c}}}(x)=\frac{\pi}{2}
-\frac{x^2}{4 \sqrt{3}} - \frac{x^4}{120 \sqrt{3}}
+\cO(x^6),
\]
and then \eqref{eq:Phi} is written as
\[
\Phi_{t_{\mathrm{c}}}(x)
= 2\arcsinh\left[\bbmi - \frac{1}{72}(\sqrt{3}+\bbmi) x^4+\cO(x^6)\right].
\]
Therefore, we can verify that 
$1-R_{t_{\mathrm{c}}}(x) \sim x^2/(3 \sqrt{2})$,
$\varphi_{t_{\mathrm{c}}}^{\mathrm{c}}-\varphi_{t_{\mathrm{c}}}(x)
\sim x^2/\sqrt{6}$, and hence
\eqref{eq:edgeB2} is proved. 
\end{proof}
%%%%%%%%%%%%%%%%%%%%%%%

Similar asymptotic analysis is possible for the edges of $\wtilde{\gamma}_{t}$ 
with $0< t\leq t_{\mathrm{c}}$.
In fact, \eqref{eq:Ftb} admits the expansion
\[
F_t(r)=1-2t+\frac{1}{6} (3-2t) (1-r)^2+
\cO((1-r)^3)
\quad \text{as}\quad r \uparrow 1.
\]
Note that the curve $\wtilde{\gamma}_{t}$ can be parametrized as 
$\wtilde{\gamma}_{t}=\{r_{t}(\theta)e^{\bbmi \theta}|\theta\in [-\theta_{t}^{\mathrm{c}},\theta_{t}^{\mathrm{c}}]\}$.
We can verify the following. 
%%%%%%%%%%%%%%%%%%%%%
\begin{prop}
\label{thm:vicinity}
The curve $\widetilde{\gamma}_t$, $0< t \leq t_{\mathrm{c}}$ behaves
as follows in the vicinity of the edge $e^{\bbmi\theta_{t}^{\mathrm{c}}}$. 
\begin{enumerate}
\item
For $0 < t < t_{\mathrm{c}}$,
\[
1- r_t(\theta) \sim a_t (\theta_t^{\mathrm{c}}-\theta)^{\nu^{\prime}}
\quad \mbox{as $\theta \uparrow \theta_t^{\mathrm{c}}$ \, with}
\quad \nu^{\prime}= \frac{1}{2}, 
\]
where the coefficient 
\[
a_t=\frac{2(t(1-t))^{1/4}}{(1-2t/3)^{1/2}}.
\]
vanishes as
\[
a_t \sim 2\sqrt{3} (1-t)^{\sigma'} \downarrow 0 \quad
\mbox{as $t \uparrow t_{\mathrm{c}}$ \, with} \quad
\sigma'=\frac{1}{4}.
\]
\item
At $t=t_{\mathrm{c}}$,
\[
1- r_{t_{\mathrm{c}}}(\theta) \sim \sqrt{3} (\theta_{t_{\mathrm{c}}}^{\mathrm{c}}-\theta)
\quad \mbox{as $\theta \uparrow 
\theta_{t_{\mathrm{c}}}^{\mathrm{c}}=\pi$}.
\]
\end{enumerate}
\end{prop}
%%%%%%%%%%%%%%%%%

By comparing \eqref{eq:varphi_t_c} and \eqref{eq:theta_c}, we always have $\varphi_t^{\mathrm{c}} > \theta_t^{\mathrm{c}}$ for $0 < t < t_{\mathrm{c}}$.
This accompanies the fact that $\nu'=1/2<3/2=\nu$.
In fact, it can be visually seen that
the edge behavior of $\gamma_t$ (Figure~\ref{fig:SLE} (a))
is sharper than that of $\widetilde{\gamma}_t$ (Figure~\ref{fig:SLE_pre} (a)).

We recall that the edge exponent $\nu=3/2$ for $0 < t < t_{\mathrm{c}}$ in \eqref{eq:edgeB1}
is the same value with the edge exponent obtained by~\cite{HottaKatori2018} in the chordal setting.
In this sense, we could say that $\nu$ is a {\it universal} exponent.

As $t$ approaches $t_{\mathrm{c}}$, the coefficient $b_t$ in the edge asymptotics for $\gamma_{t}$ ({\it susceptibility}, let us call it) diverges 
and the exponent eventually turns to unity at $t=t_{\mathrm{c}}$.
For $\wtilde{\gamma_{t}}$, the opposite is observed; the coefficient $a_{t}$ vanishes,
but the exponent at $t=t_{\mathrm{c}}$ is again unity.
Because of these characteristics, we would like to regard the singularity at $t=t_{\mathrm{c}}$ as a {\it critical phenomenon}.

Let us close this section by stressing a few aspects that we think are significant.
First, the critical phenomenon is a new feature of the radial setting;
there is no analogue in the chordal setting~\cite{HottaKatori2018}.
Second, the critical phenomenon can only be observed in the hydrodynamic limit $N\to\infty$. This is natural since a critical phenomenon occurs in an infinite system.
Lastly, even though a singularity at $t=t_{\mathrm{c}}$ is already observed for the hydrodynamic limit of circular Dyson Brownian motions when $\supp\mu_{t}$ changes its topology,
the critical phenomenon only occurs at the level of multiple SLE.
With that, we say that the singularity for $\supp\mu_{t}$ is a shadow of a critical phenomenon.

\subsection{Hydrodynamic limit from GFF perspective}
Since the circular Dyson Brownian motions as driving processes were derived from coupling with GFF, it is natural to ask about the hydrodynamic limit of GFF.
As we saw in Proposition~\ref{thm:Burgers1} and Proposition~\ref{thm:hydro_SLE},
the hydrodynamic limit of the circular Dyson Brownian motions and the limit Loewner chain $(\mSLEgAt{\infty}{t}:t\geq 0)$ are both deterministic. 
The corresponding limit of the stopping time $\tau_{U}$ with an open $U\subset \bD$ is also deterministic.

Assuming the equivalent conditions in Theorem~\ref{thm:coupling_multi_curve}, the following limit in distribution exists:
\begin{align*}
	\lim_{N\to\infty}\frac{1}{N}\,\gffB_{U,t/N} 
	&\,= \lim_{N\to\infty}\frac{1}{N}\harmMAt{N}{t/N} \\
	&\,= \frac{1}{\sqrt{\kappa}}\left(-2\int_{\bS}\arg (\mSLEgAt{\infty}{t}(\cdot)- x)\mu_{t}(dx) + \arg \mSLEgAt{\infty}{t}(\cdot)\right) \\
	&\,\eqqcolon \frac{1}{\sqrt{\kappa}}\harmHydroAt{t},\quad t\in [0,\tau_{U}],
\end{align*}
where $(\mu_{t}:t\geq 0)$ is the hydrodynamic limit of the circular Dyson Brownian motions in Proposition~\ref{thm:Burgers1}.
Note that $\frac{1}{N}\sum_{i=1}^{N}\Theta^{(i)}_{t} \to 0$ in distribution as $N\to \infty$.
The overall $1/\sqrt{\kappa}$ is factored out after the limit; this is consistent with the fact that the hydrodynamic limit does not depend on the initial choice of $\beta = 8/\kappa$.
Note that $\harmHydroAt{t}$ can be obtained as $\lim_{N\to\infty}N^{-1}\harmClMAt{N}{t/N}$ as well.
Similarly to the classical limit discussed in Section~\ref{sect:coupling_classical}, $(\harmHydroAt{t}(z):t\in [0,\tau_{U}])_{z\in U}$ is a family of integrals of motion for the Loewner chain $(\mSLEgAt{\infty}{t}:t\geq 0)$.

\bibliographystyle{alpha}
\bibliography{sle_gff}

\end{document}